\documentclass[12pt]{amsart}

\usepackage{amsmath}
\usepackage{amsbsy}
\usepackage{amssymb}
\usepackage{euscript}
\usepackage{amscd}
\usepackage{eq}

\newcommand{\bbc}{{\mathbb C}}

\newcommand{\bbq}{{\mathbb Q}}
\newcommand{\bbr}{{\mathbb R}}

\newcommand{\bbz}{{\mathbb Z}}

\newcommand{\al}{{\alpha}}

\newcommand{\be}{{\beta}}

\newcommand{\gam}{{\gamma}}

\newcommand{\del}{{\delta}}

\newcommand{\Del}{{\Delta}}

\newcommand{\vep}{{\varepsilon}}

\newcommand{\ka}{{\kappa}}

\newcommand{\lam}{{\lambda}}

\newcommand{\sig}{{\sigma}}

\newcommand{\g }{{\mathfrak g}}

\newcommand{\gi}{{\mathfrak i}}

\newcommand{\gn}{{\mathfrak n}}

\newcommand{\gp}{{\mathfrak p}}

\newcommand{\gX}{{\mathfrak X}}

\newcommand{\co}{{\mathcal O}}

\newcommand{\xspan}{{\operatorname{Span}}}
\newcommand{\aff}{{\operatorname {Aff}}}

\newcommand{\tr}{{\operatorname {Tr}}}

\newcommand{\lev}{{\operatorname{lev}}}
\newcommand{\n}{{\operatorname {N}}}

\newcommand{\gal}{{\operatorname{Gal}}}
\newcommand{\gl}{{\operatorname{GL}}}

\newcommand{\sst}{{\text{ss}}}

\newcommand{\vol}{{\operatorname{vol}}}
\newcommand{\ord}{{\operatorname{ord}}}

\newcommand{\pv}{prehomogeneous vector space}

\newcommand{\Z}{\bbz}
\newcommand{\Q}{\bbq}
\newcommand{\R}{\bbr}
\newcommand{\C}{\bbc}

\newcommand{\mk}{k^{\times}}

\newcommand{\md}{d^{\times}}

\newcommand{\akiadd}
{Department of Mathematics\\ Oklahoma State University \\
Stillwater OK 74078}
\newcommand{\kableadd}%
{Department of Mathematics\\ Cornell University\\
Ithaca NY 14853}
\newcommand{\akiemail}{yukie{@}math.okstate.edu}

\newcommand{\kableemail}{kable{@}math.cornell.edu}

\newcommand{\bmid}{\;\vrule\;}
\newcommand{\sub}{\subseteq}
\newcommand{\bk}{\backslash}
\newcommand{\ti}{\widetilde}

\newcommand{\beq}{\begin{equation}
\begin{aligned}}
\newcommand{\eeq}{\end{aligned}
\end{equation}}

\newtheorem{thm}{Theorem}[section]
\newtheorem{lem}[thm]{Lemma}
\newtheorem{cor}[thm]{Corollary}
\newtheorem{prop}[thm]{Proposition}

\theoremstyle{definition}

\newtheorem{defn}[thm]{Definition}

\theoremstyle{remark}

\newcommand{\kt}{{\widetilde{k}}}

\newcommand{\calo}{{\mathcal{O}}}

\newcommand{\nn}{\gn}

\renewcommand{\k}{\kappa}

\begin{document}

\title[Mean value theorem]
{The mean value of the product of class numbers of paired quadratic fields II}
\author{Anthony C. Kable}
\address{\kableadd}
\email{\kableemail}
\author{Akihiko Yukie}
\address{\akiadd}
\email{\akiemail}
\date{\today}
\begin{abstract}
This is the second part of a two part paper. In this part, 
we evaluate the previously unevaluated local densities 
at dyadic places which appear in the density
theorem stated in the first part.
For this purpose we introduce an invariant, the level, attached to a
pair of ramified quadratic extensions of a dyadic local field. This
invariant measures how close the fields are in their arithmetic
properties and may be of interest independent of its application here.
\end{abstract}
\maketitle

\section{Introduction} \label{intro2}

We first recall the main result of part I \cite{kable-yukie-pbh-I}
and this paper.  
If $k$ is a number field, let $\Del_k$, $h_k$, 
and $R_k$ be the absolute discriminant
(which is an integer), the class number, and the regulator,
respectively.  
We fix a number field $k$ and a 
quadratic extension $\ti k$ of $k$.
If $F\not=\ti k$ 
is another quadratic extension of $k$, let $\ti F$ be the
compositum of $F$ and $\ti k$.  Then $\ti F$ is a biquadratic 
extension of $k$ and so contains precisely three quadratic 
extensions, $\ti k$, $F$ and say, $F^*$, of $k$.  We say that 
$F$ and $F^*$ are \emph{paired}.  

For simplicity we specialize to the case $k=\Q$.
Let $\ti k = \Q(\sqrt{d_0})$ where $d_0\not=1$ is a square
free integer.  Suppose $|\Del_{\Q(\sqrt{d_0})}| 
= \prod_p p^{\ti\del_p(d_0)}$ is the prime decomposition.    
For any prime number $p$, we put 
\begin{equation*} 
E'_p(d_0) = \begin{cases} 
1-3p^{-3}+2p^{-4}+p^{-5}-2p^{-6} 
& \text{if $p$ is split in}\; \ti k, \\
(1+p^{-2})(1-p^{-2}-p^{-3}+p^{-4}) 
& \text{if $p$ is inert in}\; \ti k, \\
(1-p^{-1})(1+p^{-2}-p^{-3}+
p^{-2\ti\del_p(d_0)-2\lfloor {\ti\del_p(d_0)}/2 \rfloor -1})
& \text{if $p$ is ramified in}\; \ti k,
\end{cases}
\end{equation*} 
where $\lfloor {\ti\del_p(d_0)}/2 \rfloor$ is the largest integer
less than or equal to ${\ti\del_p(d_0)}/2$.   
We define 
\begin{equation*} 
\begin{aligned}
c_+ (d_0) & = \begin{cases} 16 & d_0>0, \\ 
               8\pi   & d_0<0, \end{cases} \quad 
c_- (d_0) = \begin{cases} 4\pi^2 & d_0>0, \\ 
               8\pi   & d_0<0, \end{cases} \\
M(d_0) & = |\Del_{\Q(\sqrt{d_0})}|^{\frac 12} 
\zeta_{\Q(\sqrt{d_0})}(2)\prod_p E'_p(d_0)
.\end{aligned}
\end{equation*} 
The following two theorems are the main results of part I and this
paper.
\begin{thm}\label{simple-mainthm1} With either choice of sign we have
\begin{equation*} 
\lim_{X\to\infty} X^{-2} 
\sum_{\substack{[F:\Q]=2, \\ 0< \pm \Del_F < X}}
h_FR_Fh_{F^*}R_{F^*} 
= c_{\pm}(d_0)^{-1} M(d_0)
.\end{equation*} 
\end{thm} 
\begin{thm}\label{simple-mainthm2} With either choice of sign we have
\begin{equation*} 
\lim_{X\to\infty} X^{-2} 
\sum_{\substack{[F:\Q]=2, \\ 0< \pm \Del_F < X}} 
h_{F(\sqrt{d_0})}R_{F(\sqrt{d_0})}
= c_{\pm}(d_0)^{-1} h_{\Q(\sqrt{d_0})} R_{\Q(\sqrt{d_0})} M(d_0)
.\end{equation*} 
\end{thm} 

For a general introduction to this problem, the reader should
see the introduction to part I.  
Our method of deriving density theorems such as 
Theorems \ref{simple-mainthm1} and \ref{simple-mainthm2}
from information on the zeta functions of \pv s is 
called the filtering process.  The filtering process
for this case was discussed in the introduction and 
sections 6 and 7 of part I.  The remaining task for us to 
finish the filtering process is to find the previously unevaluated
local densities at the dyadic places of $k$ and this is the 
main purpose of this part.  

Let $W$ be the space of binary Hermitian forms.   Our approach 
to the above theorems is based on a consideration of 
the zeta function for the following \pv{}:
\begin{equation}  \label{pbh}
G = \gl(2)_{\ti k}\times \gl(2),\quad 
V = W\otimes \aff^2
,\end{equation} 
where $\gl(2)_{\ti k}$ is regarded as a group over $k$ by
restriction of scalars and $\aff^2$ is affine $2$-space
regarded as a variety over $k$.  There is a 
relative invariant polynomial $P(x)$ of degree four
(given immediately after (3.5) in part I) and 
we put $V^{\sst}=\{x\in V\mid P(x)\not=0\}$.  

Let $v$ be a finite place of $k$, $k_v$ be the completion of 
$k$ at this place and  $K_v\sub G_{k_v}$ be the 
standard maximal compact subgroup of $G_{k_v}$.  
We assume that $\ti k_v = \ti k\otimes_{k} k_v$ is a field.
It is proved in \cite{kayu}, p. 324 that the orbit space
$G_{k_v}\bk V^{\sst}_{k_v}$ corresponds bijectively 
with the set of extensions of $k_v$ of degree one or two.  
For $x\in V^{\sst}_{k_v}$ we denote the field corresponding 
to $x$ by $k_v(x)$ and the identity component of the stabilizer of $x$
by $G_x^{\circ}$.  

In part I we selected standard representatives for the orbits in
$G_{k_v}\backslash V_{k_v}^{\sst}$ and introduced an equivalence
relation $\asymp$ on $V_{k_v}^{\sst}$ whose equivalence classes are
unions of $G_{k_v}$-orbits. These definitions will be reviewed,
respectively, in section \ref{review} and at the end of section
\ref{level}. On $V_{k_v}$ we use the additive Haar measure under which
$\vol(V_{\calo_v})=1$ and on $G_{x\,k_v}^{\circ}$ the Haar measure
described in \cite{kable-yukie-pbh-I}, Definition 5.13. We shall not
have to recall this latter definition here; all the information we
require about it will be presented at the beginning of section
\ref{stab-volume}. If $x$ is the standard representative for an orbit
in $G_{k_v}\backslash V_{k_v}^{\sst}$ then we define
\begin{equation*}
\vep_v(x)=\vol(G_{x\,k_v}^{\circ}\cap K_v)\vol(K_vx)
\end{equation*}
and
\begin{equation*}
\bar\vep_v(x)=\sum_{y\asymp x}\vep_v(y)
\end{equation*}
where the sum is over standard representatives for orbits in the
equivalence class of $x$. The \emph{local density} at $v$ is then
\begin{equation*}
E_v=\sum_{x}\vep_v(x)=\sum_{x}\bar\vep_v(x)
\end{equation*}
where the first sum is over all standard representatives for orbits in
$G_{k_v}\backslash V_{k_v}^{\sst}$ and the second over a set
containing one standard representative for an orbit in each class in
$G_{k_v}\backslash V_{k_v}^{\sst}/\asymp$. The values of
$\bar\vep_v(x)$ calculated in this paper are summarized in Tables
\ref{table-dyadic-ungrouped} and \ref{table-dyadic-grouped}. The
remaining notation used in these tables is defined in section
\ref{review} and at the end of section \ref{level}.  
The values of $\bar\vep_v(x)$ in Tables \ref{table-dyadic-ungrouped}
and \ref{table-dyadic-grouped} are verified in  
Propositions \ref{stabvolume}, \ref{orbit-volume}, \ref{orbit-volume*},
Corollary \ref{orbit-volume-ur} and \cite{kable-yukie-pbh-I}, 
Proposition 10.3.  

\newcommand{\tallstrut}{\rule[-7pt]{0cm}{22pt}}
\newcommand{\tinsert}{\quad\tallstrut}
\begin{table}
\begin{center}
\begin{tabular}{|c|c|}
\hline
\ Index\ & \tinsert$\bar\vep_v(x)$\tinsert \\
\hline\hline
(rm rm)* & \tinsert 
$\tfrac12q_v^{-2\ti\del_v-{2\lfloor {\ti\del_v}/2 \rfloor}}
(1-q_v^{-2})^2$
\tinsert \\
\hline
(rm rm ur)  &\tinsert 
$q_v^{-2\ti\del_v}
(1-\tfrac12q_v^{-{2\lfloor {\ti\del_v}/2 \rfloor}})
(1-q_v^{-1})^2(1-q_v^{-2})$
\tinsert \\
\hline
\end{tabular}
\end{center}
\vspace*{5pt}
\caption{$\bar\vep_v(x)=\vep_v(x)$ for types (rm rm)* and (rm rm ur)}
\label{table-dyadic-ungrouped}
\end{table}

\begin{table}
\begin{center}
\begin{tabular}{|c|c|}
\hline
\ Conditions\ &\tinsert$\bar\vep_v(x)$\tinsert \\
\hline\hline
$\del_{x,v}\neq\ti\del_v$, $\del_{x,v}\leq 2m_v$&\tinsert 
$q_v^{-(\del_{x,v}/2+\lam_{x,v})}
(1-q_v^{-1})^2(1-q_v^{-2})^2$
\tinsert \\
\hline
$\del_{x,v}\neq\ti\del_v$, $\del_{x,v}=2m_v+1$&\tinsert
$q_v^{-(m_v+\lam_{x,v}+1)}(1-q_v^{-1})(1-q_v^{-2})^2$
\tinsert \\
\hline
$\del_{x,v}=\ti\del_v\leq 2m_v$,
$\lam_{x,v}=\tfrac12\ti\del_v$ & \tinsert 
$q_v^{-2\lam_{x,v}}(1-q_v^{-1})(1-2q_v^{-1})(1-q_v^{-2})^2$
\tinsert \\
\hline
$\del_{x,v}=\ti\del_v\leq 2m_v$,
$\lam_{x,v}>\tfrac12\ti\del_v$ & \tinsert 
$q_v^{-2\lam_{x,v}}(1-q_v^{-1})^2(1-q_v^{-2})^2$
\tinsert \\
\hline
$\del_{x,v}=\ti\del_v=2m_v+1$ &\tinsert 
$q_v^{-2\lam_{x,v}}(1-q_v^{-1})^2(1-q_v^{-2})^2$
\tinsert \\
\hline
\end{tabular}
\end{center}
\vspace*{5pt}
\caption{$\bar\vep_v(x)$ for grouped 
dyadic orbits of type (rm rm rm)}
\label{table-dyadic-grouped}
\end{table}

All the cases we have to deal with here involve pairs $(\ti
k_v,k_v(x))$ of ramified quadratic extensions of 
$k_v$.  Since $v$ is dyadic, they are both 
wildly ramified and this is the main difficulty 
of the situation.  The definition of $\vep_v(x)$
consists of two factors, $\vol(G^{\circ}_{x\, k_v}\cap K_v)$
and $\vol(K_v x)$.  It is the second factor which requires
grouping of orbits to compute.  So, for us to be able to 
compute $\bar\vep_v(x)$, the first factor has be 
the same for all $x$ in the same group.  This means 
that the grouping has to be coarse enough to compute
the sum of the second factors, but fine enough so that 
the first factor stays constant in every group.  
When we defined the appropriate grouping in section 7 of part I, 
we used the relative discriminants of the extensions 
$k_v(x)/k_v$ and $\ti k_v(x)/\ti k_v$, where 
$\ti k_v(x)$ is the compositum of $\ti k_v$ 
and $k_v(x)$.  However, we would like to use 
congruence conditions on the vector space $V$ to 
compute the sum of $\vol(K_v x)$ and it is not 
easy to relate the relative discriminant of 
$\ti k_v(x)/\ti k_v$ directly to congruence conditions 
on $V$.  

To surmount this difficulty,
we introduce, in section \ref{review}, the notion of the \emph{level}
of a pair $(k_1,k_2)$ of ramified quadratic extensions of $k_v$. This
number provides a measure of how close $k_1$ and $k_2$ are in their
arithmetic properties and we prove that the grouping with respect to
the level is the same as
the grouping with respect to the relative discriminants of 
$k_v(x)/k_v$ and $\ti k_v(x)/\ti k_v$.  
The definition of the level itself
involves congruence conditions and so it is relatively easy 
to relate it to congruence conditions on $V$.  
After establishing the properties of the level, it is 
fairly straightforward to carry out the computation of 
$\bar{\vep}_v(x)$.  

For the rest of this introduction we discuss 
the organization of this part.  
Throughout this part, $k$ is a fixed number field,
and $\ti k$ is a fixed quadratic extension of $k$.  
We also assume throughout that $v$ is a dyadic place of 
$k$ and $\ti k_v$ is a ramified 
quadratic extension of $k_v$.  Therefore, the content
of this part is of a purely local nature.  
Even though we basically follow the notation and definitions
in part I, a minimal review of basic notions and definitions
should help the reader, and we shall provide this in section 
\ref{review}.  In section \ref{level}, we introduce the notion
of the level of two ramified quadratic extension of a dyadic
local field and establish its fundamental properties.  
For the sake of computing 
$\bar{\vep}_v(x)$, Proposition \ref{crucial-prop}
is the crucial result.  In section \ref{stab-volume}, we compute
$\vol(G_{x\,k_v}^{\circ}\cap K_v)$ and prove that it depends only on
the level of $k_v(x)$ and $\ti k_v$.
In section \ref{bxrf}, we compute the sum of 
$\vol(K_vx)$ for each equivalence class of representatives,
using the same method as that in section 11 of part I.

\section{Review of facts from part I}\label{review}

In this section we give a minimal review 
of basic notation and definitions from part I which 
are needed in this part.  

If $X$ is a finite set then $\# X$ will denote its cardinality. The
standard symbols $\Q$, $\R$, $\C$ and $\Z$ will denote respectively
the rational, real and complex numbers and the rational integers.
If $a\in\R$ then the largest integer $z$ such that $z\leq a$ is
denoted $\lfloor a\rfloor$ and the smallest integer $z$ such that
$z\geq a$ by $\lceil a\rceil$. 
If $R$ is any ring then $R^{\times}$ is the set
of invertible elements of $R$ and if $V$ is a variety defined over
$R$ then $V_R$ denotes its $R$-points. If $G$ is an algebraic group
then $G^{\circ}$ denotes its identity component.

Throughout this paper,
$k$ is a fixed number field, $\ti k$ is a fixed quadratic
extension of $k$ and $v$ is a dyadic place of $k$ such 
that $\ti k_v=\ti k\otimes_{k} k_v$ is a ramified 
quadratic extension of $k_v$.  We denote the non-trivial element of 
$\gal(\ti k/k)$ by $\sig$.  
Let $\co_v,\ti\co_v$ be the integer rings of $k_v,\ti k_v$
and $\gp_v=(\pi_v),\ti\gp_v=(\ti\pi_v)$ be their prime ideals.  
We denote the absolute value in $k_v$ by $|\;|_v$.  
As far as  notation pertaining to number fields and local fields,
we use the same conventions as in part I: the notation for the $\ti k$
object will be derived from that of the $k$ object by
adding a tilde and, for other fields, 
by writing the field in question 
as the subscript. For example, $\calo_F$ for the ring of integers of
the field $F$. 
If $a\in k_v$ and $(a)=\gp_v^i$ then we write 
$\ord_{k_v}(a)=i$. If $\gi$ is a fractional ideal in $k_v$ and
$a-b\in\gi$ then we write $a\equiv b\;(\gi)$ or $a\equiv b\;(c)$ if $c$
generates $\gi$.

If $k_1/k_2$ is a finite extension either of local fields or of
number fields then we shall write $\Del_{k_1/k_2}$ for the relative
discriminant of the extension; it is an ideal in the ring of
integers of $k_2$. We put $\Del_{\ti k_v/k_v}=\gp_v^{\ti\del_v}$.
We shall use the notation $\tr_{k_1/k_2}$ and $\n_{k_1/k_2}$ for
the trace and the norm in the extension $k_1/k_2$.

We assume that the reader is familiar with the basic definitions
and facts concerning local fields. These may be found in \cite{weilc}.
We choose Haar measures $dx_v$ on $k_v$ and $\md t_v$  $\mk_v$
so that $\int_{\co_v}dx_v=1$ and   
$\int_{\co_v^{\times}} \md t_v = 1$.  

As in part I, we use the following notation
\begin{equation} 
a(t_1,t_2)= \pmatrix t_1 & 0\\ 0 & t_2\endpmatrix,\;
n(u) = \pmatrix 1 & 0\\ u & 1\endpmatrix\,.
\end{equation}

Let $(G,V)$ be the \pv{} (\ref{pbh}) in the introduction.  
We identify $x=(x_1,x_2)\in V$ with 
the $2\times 2$-matrix $M_x(v)=v_1x_1+v_2x_2$
of linear forms in the variables $v_1$ and $v_2$, which we collect
into the row vector $v=(v_1,v_2)$. With this identification, 
the action of $g=(g_1,g_2)\in G$ on $V$  is 
$M_{gx}(v)= g_1M_x(vg_2)\,{}^{t}g_1^{\sigma}$.  
We define $F_x(v) = -\det M_x(v)$.  
Then $F_{gx}(v) = \n_{\ti k/k}(\det g_1) F_x(v g_2)$.  
It is proved in \cite{kayu}, p. 324 that by associating $x$ with the
splitting field of
$F_x(v)$, the orbit space $G_{k_v}\bk V^{\sst}_{k_v}$
corresponds bijectively with field extensions $F/k_v$ of degree
one or two.  If $x\in V^{\sst}_{k_v}$ then we denote the
corresponding field by $k_v(x)$.  If $k_v(x)\not=k_v,\ti k_v$ then
we define $\ti k_v(x)$ to be the compositum of $\ti k_v$ and $k_v(x)$.

We use coordinate systems on $G$ and $V$ 
similar to those in part I, as follows.  
For elements $g=(g_1,g_2)\in G$ we shall write
\begin{equation}\label{gform}
g_i=\begin{pmatrix}g_{i11}&g_{i12}\\g_{i21}&g_{i22}\end{pmatrix}
\end{equation}
for $i=1,2$. For vectors $x=(x_1,x_2)\in V$ we shall
put
\begin{equation}  \label{xform}
x_i=\begin{pmatrix}x_{i0}&x_{i1}\\x_{i1}^{\sigma}&x_{i2}
\end{pmatrix}
.\end{equation}
With this coordinate system,  $F_x(v) = a_0(x)v_1^2+
a_1(x)v_1v_2+a_2(x)v_2^2$ where
\begin{equation}\label{aform}
\begin{aligned}
a_0(x) & = \n_{\ti k_v/k_v}(x_{11})-x_{10}x_{12}, \\
a_1(x) & = \tr_{\ti k_v/k_v}(x_{11}x_{21}^{\sig})
-x_{10}x_{22}-x_{12}x_{20}, \\
a_2(x) & = \n_{\ti k_v/k_v}(x_{21})-x_{20}x_{22}\,
.\end{aligned}
\end{equation}

Suppose that $p(z)=z^2+a_1z+a_2\in k[z]$ has distinct roots $\al_1$
and $\al_2$. We collect these into a set $\al=\{\al_1,\al_2\}$,
since the numbering is arbitrary. Define $w_p\in V_k$ by
\begin{equation} 
w_p = \left(\pmatrix 0 & 1 \\ 1 & a_1\endpmatrix,
\pmatrix 1 & a_1\\ a_1 & a_1^2-a_2\endpmatrix\right)\,
.\end{equation} 
Then $F_{w_p}(z,1)=p(z)$ and so we can choose a
representative of the form $w_p$ for each orbit in the orbit space
$G_{k_v}\bk V^{\sst}_{k_v}$. These are the
\emph{standard representatives}. As remarked in
\cite{kable-yukie-pbh-I}, (3.15) and what follows, if we put
\begin{equation*}
w=\left(\begin{pmatrix}1&0\\0&0\end{pmatrix},
\begin{pmatrix}0&0\\0&1\end{pmatrix}\right)
\end{equation*}
and
\begin{equation*}
h_{\al}=\begin{pmatrix}1&-1\\-\al_1&\al_2\end{pmatrix}
\end{equation*}
then we have $w_p=(h_{\al},(\al_2-\al_1)^{-1}h_{\al})w$ if $k_v(w_p)\neq
\ti k_v$ and $w_p=(h_{\al},h_{\al},(\al_2-\al_1)^{-1}h_{\al})w$ if
$k_v(w_p)=\ti k_v$. (In the latter case we are regarding $G_{k_v}$ as
being embedded in $G_{\ti k_v}$; this is explained more fully in
\cite{kable-yukie-pbh-I}, section 3.) 

We only consider $x$ such that 
$k_v(x)/k_v$ is a ramified quadratic extension.  
Since $\ti k_v/k_v$  is also ramified, by assumption, 
there are three types of orbits. By definition, the
one corresponding to $\ti k_v$ has index (rm rm)*, those corresponding
to quadratic extensions $k_v(x)/k_v$ such that $k_v(x)\neq \ti k_v$ and
$\ti k_v(x)/\ti k_v$ is unramified have 
index (rm~rm~ur) and those
corresponding to quadratic extensions $k_v(x)/k_v$ such that
$k_v(x)\neq\ti k_v$ and $\ti k_v(x)/\ti k_v$ is ramified 
have index (rm~rm~rm).  These indices are used in
Tables \ref{table-dyadic-ungrouped}, \ref{table-dyadic-grouped}.

\section{The level of paired quadratic fields}\label{level}
  
Let $k_1\not=k_2$ be ramified quadratic extensions of $k_v$, 
and $k_1\cdot k_2$ be the  compositum of $k_1$ and $k_2$.   
We introduce the notion of the level and prove its fundamental
properties in this section.  For the rest of this paper we put
$2\co_v=\gp_v^{m_v}$. 

First we need to recall some facts concerning quadratic extensions of
$k_v$. There is a unique unramified quadratic extension of $k_v$ and
it is well-known that it is generated by a root of the Artin-Schreier
polynomial $z^2-z-c$ for a suitable choice of $c\in\calo_v^{\times}$.
Thus it is also generated by the square-root of $1+4c$. If
$\vep\in\calo_v^{\times}$ is a unit whose square-root generates the
unramified quadratic extension of $k_v$ then $\vep=a^2(1+4c)$ for some
$a\in\calo_v^{\times}$ and so the congruence $\vep\equiv
a^2\;(\gp_v^{2m_v})$ is solvable. Conversely, if
$\vep\in\calo_v^{\times}$ is such that $\vep\equiv
a^2\;(\gp_v^{2m_v})$ is solvable but $\vep\equiv
a^2\;(\gp_v^{2m_v+1})$ is not, then $\vep$ is not a square and
$(2a)^{-1}(a-\sqrt{\vep})$ is easily seen to satisfy an Artin-Schreier
polynomial, so that $\sqrt{\vep}$ generates the unramified quadratic
extension of $k_v$. Notice that $\vep\equiv a^2\;(\gp_v^{2m_v+1})$
being solvable implies that $\vep$ is a square, by Hensel's lemma.

Now we turn to ramified quadratic extensions, $F$, of $k_v$. Every
such extension is generated by a root of an Eisenstein polynomial
$p(z)=z^2+a_1z+a_2$. This root is a uniformizer, $\pi_F$, of $F$ and
we have $\calo_F=\calo_v[\pi_F]$ and hence
$\Del_{F/k_v}=(a_1^2-4a_2)\calo_v$ for any choice of Eisenstein
polynomial which splits in $F$. If $\ord_{k_v}(a_1)\geq m_v+1$ then we
may make the transformation $z\mapsto z-(a_1/2)$ in order to assume
that $a_1=0$. These extensions are exactly those generated by the
square-root of a uniformizer of $k_v$ and they have
$\Del_{F/k_v}=\gp_v^{2m_v+1}$. If $1\leq\ord_{k_v}(a_1)\leq m_v$ then
put $\ell=\ord_{k_v}(a_1)$. Here $\Del_{F/k_v}=\gp_v^{2\ell}$ and $F$
is generated by the square-root of $a_1^2-4a_2$ and hence also by the
square-root of the unit $1-4a_2a_1^{-2}=1+\pi_v^{2(m_v-\ell)+1}c$ for
a suitable $c\in\calo_v^{\times}$. This exhausts all quadratic
extensions of $k_v$. If $\vep\in\calo_v^{\times}$ is a non-square unit
and $\vep\equiv a^2\;(\gp_v^{2m_v})$ is not solvable then let $i<2m_v$
be the largest integer such that $\vep\equiv a^2\;(\gp_v^{i})$ is
solvable. We must have $\vep = a^2(1+\pi_v^{2(m_v-\ell)+1}c)$ for some
$1\leq\ell\leq m_v$ and $c\in\calo_v^{\times}$ and then
$i=2(m_v-\ell)+1$. In this case, $\pi_v^{\ell-m_v}(\sqrt{\vep}-a)$ is
a uniformizer of $k_v(\sqrt{\vep})$. From this paragraph and the
previous one it follows that if $\vep\in\calo_v^{\times}$ is a
non-square unit then we may always multiply $\vep$ by a square to
arrange either $\vep=1+4c$ or $\vep=1+\pi_v^{2(m_v-\ell)+1}c$ with
$c\in\calo_v^{\times}$.

In what follows we shall use the subscript $1$ (resp. $2$) to denote
objects associated with $k_1$ (resp. $k_2$). Thus $\calo_1$ will be
the ring of integers of $k_1$, $\pi_1$ a uniformizer of $k_1$, $\gp_1$
the prime ideal in $\calo_1$ and $\Del_{k_1/k_v}=\gp_v^{\del_1}$ and
similarly with $1$ replaced by $2$. Let $p_1(z)=z^2+a_1z+a_2$ and
$p_2(z)=z^2+b_1z+b_2$ be the minimal polynomials of $\pi_1$ and
$\pi_2$ over $k_v$, respectively. Let $\ell_1=\ord_{k_v}(a_1)$ if this
is less than or equal to $m_v$ and $\ell_1=m_v+1$ otherwise. Define
$\ell_2$ similarly for $k_2$. Notice that we have $\ell_i=\lfloor
(\del_i+1)/2\rfloor$.

In the following two lemmas, $F/k_v$ is a ramified quadratic
extension, $\gp_F$ is the maximal ideal in the ring of integers of $F$
and $\Del_{F/k_v}=\gp_v^{\del_F}$. We let $\ell_F=\lfloor
(\del_F+1)/2\rfloor$.

\begin{lem}\label{tracelemma}
Suppose $x\in F$ and $\ord_{F}(x)=1$. 
Then $\tr_{F/k_v}(x) \in \gp_v^{\ell_F}$.  
Moreover, if $\ell_F\leq m_v$ then
$\ord_{k_v}(\tr_{F/k_v}(x))=\ell_F$.  
\end{lem} 
\begin{proof} 
We have $\calo_F=\calo_v[x]$ and so if $z^2+c_1z+c_2$ is the minimal
polynomial of $x$ over $k_v$ then $c_1=-\tr_{F/k_v}(x)$ and
$(c_1^2-4c_2)\calo_v=\Del_{F/k_v}$. If $\ell_F\leq m_v$ then
$\Del_{F/k_v}=\gp_v^{2\ell_F}$ and hence $\ord_{k_v}(c_1)=\ell_F$. If
$\ell_F=m_v+1$ then $\Del_{F/k_v}=\gp_v^{2m_v+1}$ and so
$c_1^2\in\gp_v^{2m_v+1}$, which gives $c_1\in\gp_v^{\ell_F}$.
\end{proof}
\begin{lem}\label{jorder}
Suppose $u\in F$ and $\ord_{F}(u) = j$.  Then  
\begin{equation*}
\ord_{k_v}(\tr_{F/k_v}(u))\geq\lfloor(j+\del_F)/2\rfloor\,.
\end{equation*}
\end{lem}
\begin{proof} 
The different of $F/k_v$ is $\gp_F^{\del_F}$ and so, from the
definition of the different, $u\in\gp_F^{-\del_F}$ implies that
$\tr_{F/k_v}(u)\in\calo_v$. Multiplying by $\pi_v^n$, we find that
$u\in\gp_F^{2n-\del_F}$ implies that $\tr_{F/k_v}(u)\in\gp_v^n$. Let
$n=\lfloor(j+\del_F)/2\rfloor$. Then $2n\leq j+\del_F$ and so
$2n-\del_F\leq j$. Thus $u\in\gp_F^{2n-\del_F}$ and so
$\tr_{F/k_v}(u)\in\gp_v^n$.
\end{proof} 

For $0\leq i_1\leq i_2\leq i_1+1$ we define
\begin{equation} 
S_{i_1,i_2}(k_1,k_2) 
= \left\{\eta \in  \co_2/ \pi_2^{i_1+i_2}\co_2 \bmid
\begin{matrix} \tr_{k_2/k_v}(\eta) \equiv a_1 \; (\gp_v^{i_1}), \\
\n_{k_2/k_v}(\eta) \equiv a_2 \; (\gp_v^{i_2}) \end{matrix} \right\}
.\end{equation} 

We first show that the conditions defining $S_{i_1,i_2}(k_1,k_2)$
depend only on the class of $\eta$ modulo $\pi_2^{i_1+i_2}\calo_2$, so
that the definition makes sense. Suppose that $\eta\in\calo_2$ and
$u\in\pi_2^{i_1+i_2}\calo_2$. Then, by Lemma \ref{jorder}, 
\begin{equation*}
\ord_{k_v}(\tr_{k_2/k_v}(u))\geq\lfloor(i_1+i_2+\del_2)/2\rfloor
\geq\lfloor(2i_1+\del_2)/2\rfloor\geq i_1
\end{equation*}
and so $\tr_{k_2/k_v}(\eta)\equiv\tr_{k_2/k_v}(\eta+u)\;
(\gp_v^{i_1})$. Also,
\begin{equation*}
\n_{k_2/k_v}(\eta+u)=\n_{k_2/k_v}(\eta)+\tr_{k_2/k_v}(\eta^{\sigma}u)+
\n_{k_2/k_v}(u)
\end{equation*}
and
\begin{equation*}
\ord_{k_v}(\tr_{k_2/k_v}(\eta^{\sigma}u))\geq
\lfloor(i_1+i_2+\del_2)/2\rfloor\geq
\lfloor(2i_2+\del_2-1)/2\rfloor\geq i_2
\end{equation*}
by Lemma \ref{jorder} and the fact that $\del_2\geq2$. Further,
$\ord_{k_v}(\n_{k_2/k_v}(u))=i_1+i_2\geq i_2$ and so
$\n_{k_2/k_v}(\eta+u)\equiv\n_{k_2/k_v}(\eta)\;(\gp_v^{i_2})$. We
shall, by a slight abuse of notation, confuse elements of $\calo_2$
with their classes modulo $\pi_2^{i_1+i_2}\calo_2$, so that we may
write $\eta\in S_{i_1,i_2}(k_1,k_2)$ if the class of $\eta\in\calo_2$
satisfies the indicated conditions.

We let $\nn_1(k_1,k_2,i)$ (resp. $\nn_2(k_1,k_2,i)$) be the
cardinality of the set $S_{i,i}(k_1,k_2)$ (resp. $S_{i,i+1}(k_1,k_2)$)
for $i\geq0$. The set $S_{i_1,i_2}(k_1,k_2)$ depends on the choice of
an Eisenstein polynomial for $k_1$. However, it is $\nn_1(k_1,k_2,i)$
and $\nn_2(k_1,k_2,i)$ which interest us and it turns out that these
numbers depend only on $k_1$, $k_2$ and $i$, as we show in Lemma
\ref{independent} below. In fact, we are really only interested in the
range of $i$ in which $\nn_1(k_1,k_2,i)$ and $\nn_2(k_1,k_2,i)$ do not
vanish and Lemma \ref{independent} is more than we require. We shall
discuss the motivation for our approach at the end of this section.

\begin{defn} The largest integer $i$ such that
$S_{i,i}(k_1,k_2)\neq\emptyset$ will be called the level of $k_1$ and
$k_2$ and denoted by $\lev(k_1,k_2)$. 
\end{defn}

Of course, $\lev(k_1,k_2)$ is the largest integer, $i$, such that
$\nn_1(k_1,k_2,i)\neq0$. It is an easy consequence of Hensel's lemma
that $\lev(k_1,k_2)<\infty$ since $k_1$ and $k_2$ are distinct. A
specific upper bound for $\lev(k_1,k_2)$ will be given in Proposition
\ref{idef}. 
It follows directly from the definition that
\begin{equation}\label{easycases}
\nn_1(k_1,k_2,0)=\nn_2(k_1,k_2,0)=1,\; 
\nn_1(k_1,k_2,1)=q_v
.\end{equation}
\begin{lem}\label{independent}
\begin{itemize}
\item[(1)]
The numbers $\nn_1(k_1,k_2,i)$ and $\nn_2(k_1,k_2,i)$ depend only on
$k_1$ and $k_2$, not on the particular choice of Eisenstein
polynomial used to evaluate them. Thus 
this notation is legitimate. 
\item[(2)] For $j=1,2$, we have $\nn_j(k_2,k_1,i)= 
\nn_j(k_1,k_2,i)$ for all $i\geq0$. 
\end{itemize}
\end{lem}
\begin{proof}
If $\pi_1$ and $\pi_1'$ are uniformizers of $k_1$
then $\pi_1=c+d\pi_1'$
with $c\in \gp_v$ and $d\in\co_v^{\times}$. If
$p_1(z)=z^2+a_1z+a_2$ is the 
Eisenstein polynomial associated to
$\pi_1$ then the Eisenstein polynomial,
$p_1'(z)=z^2+a_1'z+a_2'$, associated to $\pi_1'$ is
$p_1'(z)=z^2+d^{-1}(a_1+2c)z+d^{-2}
(c^2+a_1c+a_2)$. Say $\eta\in \calo_2$ satisfies the
congruences $\tr_{k_2/k_v}(\eta) \equiv a_1\; (\gp_v^{i_1})$ and
$\n_{k_2/k_v}(\eta)\equiv a_2\; (\gp_v^{i_2})$. 
Then it is easy to check
that $\eta'=d^{-1}(\eta+c)$ satisfies the congruences
$\tr_{k_2/k_v}(\eta')\equiv a_1'\; (\gp_v^{i_1})$ 
and $\n_{k_2/k_v}(\eta')\equiv a_2'\; (\gp_v^{i_2})$. 
Since $d\in\calo_v^{\times}$, the map $\eta\mapsto d^{-1}(\eta+c)$
induces a well-defined map on $\calo_2/\pi_2^{i_1+i_2}\calo_2$ with
inverse induced by $\eta'\mapsto d\eta'-c$. This establishes a
one-to-one correspondence between the two sets and (1) follows.

Fix a uniformizer $\pi_1$ of $k_1$ and let $p_1(z)
= z^2+a_1z+a_2$ be the corresponding Eisenstein 
polynomial.  Consider $S_{i_1,i_2}(k_1,k_2)$.  
We may assume that $i_2\geq 2$, since we have
evaluated the numbers $\nn_1(k_1,k_2,0)$, $\nn_1(k_1,k_2,1)$ and
$\nn_2(k_1,k_2,0)$ in (\ref{easycases})
and they satisfy the second claim. With this assumption, every element
of $S_{i_1,i_2}(k_1,k_2)$ is (the class of) a uniformizer in
$\calo_2$.

Suppose $S_{i_1,i_2}(k_1,k_2)\not=\emptyset$.  
Fix $\eta_0\in S_{i_1,i_2}(k_1,k_2)$. 
We will use the corresponding Eisenstein
polynomial $p_0(z)=z^2+a_{01}z+a_{02}$ to evaluate $\nn_j(k_2,k_1,i)$.
Every other element $\eta$ of $S_{i_1,i_2}(k_1,k_2)$ has the form
$\eta=c(\eta)+d(\eta)\eta_0$ with
$c(\eta)\in\gp_v$ and $d(\eta)\in\co_v^{\times}$.
Moreover, the conditions on $\eta$ imply that $c(\eta)$ and $d(\eta)$
satisfy the congruences
\begin{align*}
-d(\eta)(a_{01}-2c(\eta)d(\eta)^{-1})
&\equiv a_1\; (\gp_v^{i_1}) \\
c(\eta)^2-a_{01}c(\eta)d(\eta)+a_{02}d(\eta)^2
&\equiv a_2\; (\gp_v^{i_2})\,.
\end{align*}
We define
$\varpi(\eta)=d(\eta)^{-1}(\pi_1+c(\eta))$.  Then, using the 
facts that $c(\eta)\in \gp_v$ and $i_2\leq i_1+1$, 
it is easy to check that 
\begin{equation*}
\operatorname{Tr}_{k_1/k_v}
(\varpi(\eta))\equiv a_{01}\; (\gp_v^{i_1}),\; 
\operatorname{N}_{k_1/k_v}(\varpi(\eta))\equiv
a_{02}\; (\gp_v^{i_2})
\end{equation*}
and so $\varpi(\eta)\in S_{i_1,i_2}(k_2,k_1)$.
Suppose $u\in \pi_2^{i_1+i_2}\calo_2$ and $\eta' = \eta+u$.
If we write  $u = c(u)+d(u)\eta_0$ with $d(c),d(u)\in \co_v$
then $c(u)\in \gp_v^{i_2}$ and $d(u)\in \gp_v^{i_1}$.  
By computation, 
\begin{align*}
& d(\eta')^{-1}(\pi_1+c(\eta'))- 
d(\eta)^{-1}(\pi_1+c(\eta)) \\
& = d(\eta')^{-1}d(\eta)^{-1}(-d(u)\pi_1
+c(u)d(\eta)-c(\eta)d(u))
.\end{align*}
It is easy to check that this element belongs to 
$\pi_1^{i_1+i_2}\co_1$ and so 
the map $\eta\mapsto\varpi(\eta)$ induces a well-defined map from
$S_{i_1,i_2}(k_1,k_2)$ to $S_{i_1,i_2}(k_2,k_1)$. Reversing the roles
of $k_1$ and $k_2$ we obtain a similar map from $S_{i_1,i_2}(k_2,k_1)$
to $S_{i_1,i_2}(k_1,k_2)$ induced by the map sending
$\zeta=c'(\zeta)+d'(\zeta)\pi_1$ to
$d'(\zeta)^{-1}(\eta_0+c'(\zeta))$. It is easy to check that these
maps are inverse to one another and so $S_{i_1,i_2}(k_1,k_2)$ and
$S_{i_1,i_2}(k_2,k_1)$ have the same cardinality.
\end{proof}

Let $k_3$ be the unique quadratic extension of $k_v$ other than 
$k_1$ and $k_2$ contained in $k_1\cdot k_2$.  
Let $p_1(z),p_2(z)$ be as before.  
Let $\al_1$ and $\al_2$ be the roots of $p_1$ 
and $\be_1$ and $\be_2$ be the roots of $p_2$. 

Define
\begin{align*}
\gam_1&=(\al_1-\be_1)(\al_2-\be_2), \\
\gam_2&=(\al_1-\be_2)(\al_2-\be_1)\,.
\end{align*}

The following lemma provides an equation defining $k_3$. 
We will not provide the proof since it is elementary.
\begin{lem}  \label{k3generator}
The numbers $\gam_1$ and $\gam_2$ 
generate $k_3$ over $k_v$ and are the
roots of the polynomial
\begin{equation} \label{p3defn}
p_3(z)=z^2-[2(a_2+b_2)-a_1b_1]z+J
\end{equation}
where $J=(a_2-b_2)^2+(a_1-b_1)(a_1b_2-a_2b_1)$. Moreover,
$\gam_1-\gam_2=(\al_1-\al_2)(\be_1-\be_2)$.
\end{lem}

Next we consider the relation between discriminants of $k_1,k_2,k_3$.
Let $\Del_{k_i/k_v}=\gp_v^{\del_i}$ for $i=1,2,3$.  
Note that for $i=1,2$, $\del_i=2,\dots,2m_v$ or $2m_v+1$.  
\begin{lem} \label{bqbound}
We have $\del_3\leq \operatorname{max}\{\del_1,\del_2\}$.  
Moreover, equality holds if $\del_1\not= \del_2$.
\end{lem}
\begin{proof}  
There are two cases to consider. If two of the fields are generated
by adjoining the square-root of a uniformizer then they have equal
discriminants and the third field has a smaller discriminant (since
it is obtained by adjoining the square root of a unit). 
Therefore, we have the statement of this lemma in this case.  Otherwise,
all the fields are obtained by adjoining the square root of a unit.
Let $\vep_1$, $\vep_2$ and $\vep_3$ be the units whose square roots
generate $k_1,k_2$ and $k_3$ respectively.  
We may assume that
$\vep_j=1+\pi_v^{2(m_v-\ell_j)+1}c_j$ where $c_j\in\co_v^{\times}$ and
$\del_j=2\ell_j$ for $j=1,2$.  We may also assume $\vep_3=\vep_1\vep_2$. 
Then 
\begin{equation*}
\vep_3 = 1+ \pi_v^{2(m_v-\ell_1)+1}c_1+\pi_v^{2(m_v-\ell_2)+1}c_2
+\pi_v^{2(m_v-\ell_1)+2(m_v-\ell_2)+2}c_1c_2
.\end{equation*}

If $\ell_1>\ell_2$ then 
\begin{equation*}
\vep_3\equiv \vep_1 \; (\gp_v^{2(m_v-\ell_1)+2})
.\end{equation*}
and so $\del_3 = 2\ell_1=\del_1$.  The case $\ell_2>\ell_1$ is similar.  
If $\ell_1=\ell_2$ then 
$\vep_3\equiv 1 \; (\gp_v^{2(m_v-\ell_1)+1})$.  
If $\vep_3\equiv 1\;(4)$ then $\del_3=0$ and the inequality holds
true. Otherwise, the largest number, $i$, such that
$\vep_3\equiv1\;(\gp_v^i)$ has the form $i=2(m_v-\ell_3)+1$ with
$\ell_3\leq\ell_1=\ell_2$. Then $\del_3=2\ell_3\leq\del_1=\del_2$ and
again the inequality is true.
\end{proof}
\begin{lem}  \label{level-bound}
We have 
\begin{equation*}
2\lfloor\tfrac12\operatorname{ord}_{k_v}(J)\rfloor\leq
\del_1+\del_2-\del_3\,.
\end{equation*}
\end{lem}
\begin{proof}
Let $a=\lfloor\tfrac12\operatorname{ord}_{k_v}(J)\rfloor$ so that
$J/\pi_v^{2a}$ is either a unit of a uniformizer of $k_v$. Since
$\operatorname{N}_{k_3/k_v}(\gam_j/\pi_v^a)=J/\pi_v^{2a}$ 
for $j=1$ and $2$,
we conclude that $\gam_j/\pi_v^a$ is an integer.
Thus the ideal generated by $(\gam_1-\gam_2)^2/\pi_v^{2a}$ 
in $\co_v$ is contained in $\gp_v^{\del_3}$. But
$\gam_1-\gam_2=(\al_1-\al_2)(\be_1-\be_2)$ and so the ideal generated
by $(\gam_1-\gam_2)^2$ is $\gp_v^{\del_1+\del_2}$. The inequality follows.
\end{proof}
\begin{prop}\label{idef}  
We have $2\lev(k_1,k_2)+\del_3\leq\del_1+\del_2$.
\end{prop}
\begin{proof}  Let $i=\lev(k_1,k_2)$.  
We choose Eisenstein polynomials $p_1,p_2$ 
so that $a_1\equiv b_1,\; a_2\equiv b_2 \;(\gp_v^i)$.
Then 
\begin{equation} \label{J-value}
J=(a_2-b_2)^2+(a_1-b_1)[a_1(b_2-a_2)+a_2(a_1-b_1)]
\end{equation}
and our assumptions imply that this lies in $\gp_v^{2i}$. Using the
previous lemma we obtain $2i\leq \del_1+\del_2-\del_3$
and the inequality follows. 
\end{proof}
\begin{cor}\label{ib}
\begin{itemize}
\item[(1)] If $\ell_1\not=\ell_2$ then
$\lev(k_1,k_2)\leq \operatorname{min}\{\ell_1,\ell_2\}$.
\item[(2)] If $\ell_1=\ell_2=\ell,\; \del_1=\del_2=\del$
then $\lev(k_1,k_2)\leq \del$.
\end{itemize}  
\end{cor}
\begin{proof}  Consider (1).  
Suppose, without loss of generality, that $\ell_1<\ell_2$. Then,
according to Lemma \ref{bqbound}, 
we must have $\del_3=\del_2$ and so the inequality in  
Proposition \ref{idef}  becomes $\lev(k_1,k_2)\leq \tfrac12 \del_1$. 
Since $\del_1\not= 2m_v+1$, $\tfrac12 \del_1=\ell_1$. 
Statement (2) is obvious from Proposition \ref{idef}
because $\del_3\geq0$.  
\end{proof}
Note the above corollary implies
that if $k_1,k_2/k_v$ are ramified quadratic extensions,
$\del_1=\del_2=\del$ and $S_{\del+1,\del+1}(k_1,k_2)\not=\emptyset$
then $k_1=k_2$.

\begin{prop}\label{unramified-char}
The extension $(k_1\cdot k_2)/k_2$  
is unramified if and only if $\del_1=\del_2$ and 
$S_{\del_1,\del_1}(k_1,k_2) \not=\emptyset$.  
Moreover, if these conditions are satisfied then $k_3/k_v$ is
unramified.
\end{prop} 
\begin{proof} 
Suppose $k_1 = k_v(\sqrt{\vep_1})$ and $k_2 = k_v(\sqrt{\vep_2})$.
We first assume $(k_1\cdot k_2)/k_2$ is unramified.    
Then $(k_1\cdot k_2)/k_v$ is not totally ramified.  Therefore, 
by \cite{weilc}, Corollary 4, p. 19, $k_1\cdot k_2$ 
contains an unramified 
quadratic extension of $k_v$.  Since $k_1,k_2$ are ramified, 
the remaining quadratic subfield $k_3=k_v(\sqrt{\vep_2\vep_1^{-1}})$
must be unramified over $k_v$. Let $\vep_3=\vep_2\vep_1^{-1}$, so that
$\vep_2=\vep_1\vep_3$. Multiplying $\vep_2$ and hence $\vep_3$ by a
square, if necessary, we may assume that $\vep_3\equiv1\;(4)$. Then
$\vep_1$ and $\vep_2$ have the same order in $k_v$ and, multiplying
them both by the same square, we may assume that they are either both
units or both uniformizers without altering $\vep_3$.

If $\vep_1,\vep_2$ are both uniformizers then
$\del_1=\del_2=2m_v+1$.  By the assumption on $\vep_3$,
$\vep_2=\vep_1(1+4c_3)$ for some $c_3\in\calo_v^{\times}$. 
Let $\eta_i=\sqrt{\vep_i}$ for 
$i=1,2$.  Then $\eta_1,\eta_2$ are uniformizers of 
$k_1,k_2$ respectively and 
\begin{align*}
\tr_{k_2/k_v}(\eta_2) & = \tr_{k_1/k_v}(\eta_1) = 0, \\ 
\n_{k_2/k_v}(\eta_2)  & = -\vep_2  = -\vep_1 - 4\vep_1c_3
\equiv \n_{k_1/k_v}(\eta_1) \; (\gp_v^{2m_v+1})
.\end{align*}
This implies that $S_{2m_v+1,2m_v+1}(k_1,k_2)\not=\emptyset$. 

Suppose $\vep_1,\vep_2$ are both units. Then 
$\del_1= 2\ell_1,\;\del_2=2\ell_2$ with $1\leq \ell_1,\ell_2\leq m_v$. 
Let $\vep_1= 1 + \pi_v^{2(m_v-\ell_1)+1}c_1$ and   
$\vep_3 = 1 + 4c_3$ with $c_1\in \co_v^{\times},c_3\in \co_v^{\times}$.  
Then
\begin{equation*}
\vep_2  =\vep_1\vep_3 
= 1 + \pi_v^{2(m_v-\ell_1)+1} 
(c_1+(4\pi_v^{-2m_v})\pi_v^{2\ell_1-1}c_3 +4c_1c_3)
.\end{equation*}
Let 
\begin{equation} \label{c1}
c_2 =  c_1+(4\pi_v^{-2m_v})\pi_v^{2\ell_1-1}c_3 +4c_1c_3
.\end{equation}
Then $c_2\in \co_v^{\times}$,
$\vep_2= 1 +\pi_v^{2(m_v-\ell_1)+1}c_2$ and $c_2\equiv
c_1\;(\gp_v^{2\ell_1-1})$.
Therefore, $\ell_1=\ell_2$ and so $\del_1=\del_2$.    

Let $\del=\del_1=\del_2$ and $\ell=\ell_1=\ell_2$.  
We put $\eta_i = (\pi_v^{\ell}/ 2) (\sqrt{\vep_i}-1)$ for $i=1,2$.  
Then $\eta_i$ is a uniformizer of $k_i$ satisfying 
the Eisenstein equation $z^2+\pi_v^{\ell}z-\pi_v\theta c_i=0$ for
$i=1,2$ where $\theta=\pi_v^{2m_v}/4\in\calo_v^{\times}$. Thus
\begin{align*}
&\tr_{k_2/k_v}(\eta_2)=-\pi_v^{\ell}=\tr_{k_1/k_v}(\eta_1), \\
&\n_{k_2/k_v}(\eta_2)=-\pi_v\theta c_2\equiv -\pi_v\theta
c_1\;(\gp_v^{2\ell})
\end{align*}
and since $\n_{k_1/k_v}(\eta_1)=-\pi_v\theta c_1$, we have $\eta_2\in
S_{\del,\del}(k_1,k_2)$.

Conversely, suppose $\del_1=\del_2$ and 
$S_{\del_1,\del_1}(k_1,k_2)\not= \emptyset$.  
Let $k_3,\del_3$ be as before.  
Then by Proposition \ref{idef}, $\del_3=0$.  
This implies that $k_3/k_v$ is unramified.  Since 
$k_3$ is generated by roots of an Artin-Schreier equation 
and they also generate the field extension $(k_1\cdot k_2)/k_2$, 
this extension is unramified also.  
\end{proof}
Note that by Proposition \ref{unramified-char}, 
there is precisely one orbit having (rm~rm~ur) 
as its index.  

We shall next prove that $\text{lev}(k_1,k_2)\geq 
\operatorname{min}\{
\lfloor\tfrac12(\del_1+1)\rfloor,
\lfloor\tfrac12(\del_2+1)\rfloor\}$.  
\begin{lem}\label{goup}
Suppose that $\ell_1\not=\ell_2$ and 
$1\leq i\leq \operatorname{min}\{\ell_1,\ell_2\}$
or that $1\leq i< \ell=\ell_1=\ell_2$.    
If $\eta\in \calo_2$ satisfies 
$\ord_{k_2}(\eta)=1$ and 
$\n_{k_2/k_v}(\eta)\equiv a_2 \; (\gp_v^{i})$ 
then there exists  a unit $t=c-d\eta$ 
such that $\n_{k_2/k_v}(t\eta)\equiv a_2 \; (\gp_v^{i+1})$.   
\end{lem} 
\begin{proof}  
If $i=1$, we choose $t\in\co_v^{\times}$.  
Then $\n_{k_2/k_v}(t\eta)= t^2\n_{k_2/k_v}(\eta)$. 
Since any element of $\co_v^{\times}$ is a square 
modulo $\gp_v$, we can choose $t$ so that 
$t^2\n_{k_2/k_v}(\eta)\equiv a_2\; (\gp_v^2)$. 

We now assume $i\geq 2$.  
Note that if $\ell_1\not=\ell_2$ then 
$\ell_1\leq m_v$ or $\ell_2\leq m_v$ and 
so $i+1\leq m_v+1$. This condition is obviously satisfied
in the second case.

Suppose $\eta^2+b_1'\eta+b_2'=0$ is the 
Eisenstein equation satisfied by $\eta$. 
Let  $\n_{k_2/k_v}(\eta) = b_2' = a_2+ e\pi_v^{i}$.  
Then
\begin{align*} 
\n_{k_2/k_v}(t\eta) & =\n_{k_2/k_v}(t) 
\n_{k_2/k_v}(\eta) \\
& = (c^2+ b_1'cd+b_2'd^2)(a_2+e\pi_v^{i}) \\
& \equiv (c^2 + b_2'd^2)(a_2+e\pi_v^{i})\;(\gp_v^{i+1}) \\
& \equiv a_2c^2+ a_2b_2'd^2+c^2e\pi_v^{i} \; (\gp_v^{i+1})
.\end{align*}
Note that since $i\leq \ell_2$ in both cases 
$b_1'\pi_v\equiv 0 \; (\gp_v^{i+1})$.  
Let $c= 1+\pi_v^Nf$ with $N>0$ and $f\in \co_v^{\times}$.
Then 
\begin{equation*}
c^2=1+2\pi_v^Nf+\pi_v^{2N}f^2 
\equiv 1+\pi_v^{2N}f^2 \; (\gp_v^{i+1})
.\end{equation*}
The last congruence is satisfied because of 
the condition  $i+1\leq m_v+1$.  
So 
\begin{align*}
\n_{k_2/k_v}(t\eta) & \equiv a_2+ a_2f^2\pi_v^{2N} + a_2b_2'd^2
+e\pi_v^{i}+\pi_v^{i+2N}ef^2 \\
& \equiv a_2+ a_2f^2\pi_v^{2N} + a_2b_2'd^2+e\pi_v^{i} \;(\gp_v^{i+1})
.\end{align*}
Note that the orders of $a_2f^2\pi_v^{2N}, a_2b_2'd^2$
are odd and even respectively and they can be any odd or even integer
greater than or equal to two.   
So we can choose suitable $d,f,N$ so that 
\begin{equation*}
a_2f^2\pi_v^{2N} + a_2b_2'd^2+e\pi_v^{i}
\equiv 0 \; (\gp_v^{i+1})
.\end{equation*}
\end{proof}

The following proposition provides a lower bound for 
$\lev(k_1,k_2)$.    
\begin{prop}\label{vanish}
Suppose $1\leq i\leq \operatorname{min}\{\ell_1,\ell_2\}$. 
Then $S_{i,i}(k_1,k_2)\not=\emptyset$ and so $\lev(k_1,k_2)
\geq  \operatorname{min}\{\ell_1,\ell_2\}$.  Moreover, 
if $\ell_1\not=\ell_2$ then $\lev(k_1,k_2)
=  \ell= \operatorname{min}\{\ell_1,\ell_2\}$ and 
$S_{\ell,\ell+1}(k_1,k_2)\not=\emptyset$.  
\end{prop}
\begin{proof}
We put $\ell=\min\{\ell_1,\ell_2\}$. Let $\eta\in\calo_2$ be any
uniformizer. Using Lemma \ref{goup} we can arrange that
$\n_{k_2/k_v}(\eta)\equiv a_2\;(\gp_v^{\ell})$ and
$\n_{k_2/k_v}(\eta)\equiv a_2\;(\gp_v^{\ell+1})$ if
$\ell_1\neq\ell_2$. By Lemma \ref{tracelemma},
$\tr_{k_2/k_v}(\eta)\in\gp_v^{\ell_2}\subset\gp_v^{\ell}$ and so
$\tr_{k_2/k_v}(\eta)\equiv a_1\;(\gp_v^{\ell})$. Thus $\eta\in
S_{\ell,\ell}(k_1,k_2)$ and $\eta\in S_{\ell,\ell+1}(k_1,k_2)$ if
$\ell_1\neq\ell_2$. When $\ell_1\neq \ell_2$, the equality
$\lev(\ell_1,\ell_2)=\ell$ then follows from Corollary \ref{ib}.
\end{proof}
Note that $\ell_i = \lfloor\tfrac12(\del_i+1)\rfloor$ 
for $i=1,2$ and so the above lower bound is 
$\operatorname{min}\{
\lfloor\tfrac12(\del_1+1)\rfloor,
\lfloor\tfrac12(\del_2+1)\rfloor\}$.  

\begin{lem} \label{etaform}
Suppose $i\geq 1$, $\eta,\eta'\in \co_2$,  and 
$\n_{k_2/k_v}(\eta)\equiv 
\n_{k_2/k_v}(\eta')\equiv a_2 \; (\gp_v^{i+1})$. 
Then there exist $e,f\in \co_v$ such that 
$\eta'=e\pi_v+(1+f\pi_v)\eta$.  
\end{lem}
\begin{proof}  Note that $\eta,\; \eta'$ are both uniformizers. So
we may assume that 
$\eta'=c+d\eta$ with $c\in \gp_v,d\in \co_v^{\times}$.  
Then  $\n_{k_2/k_v}(\eta') \equiv d^2 \n_{k_2/k_v}(\eta)
\equiv d^2a_2 \; (\gp_v^2)$.  Therefore, $d^2\equiv 1\; (\gp_v)$.  
This implies $d\equiv 1 \; (\gp_v)$.  
\end{proof} 

In the following proposition and its corollary
we assume that 
$\ell_1=\ell_2=\ell,\; \del_1=\del_2=\del$ and  
$\ell\leq i<\del$.  
\begin{prop} \label{crucial-prop}  
\begin{itemize}
\item[(1)] Suppose $\eta\in \co_2$ satisfies
$\tr_{k_2/k_v}(\eta)\equiv a_1\; (\gp_v^i)$,
$\n_{k_2/k_v}(\eta)\equiv a_2 \; (\gp_v^i)$. Then 
there exists $\eta'\in \co_2$ such that 
$\tr_{k_2/k_v}(\eta')\equiv a_1\; (\gp_v^i),\;
\n_{k_2/k_v}(\eta')\equiv a_2 \; (\gp_v^{i+1})$.
\item[\rm (2)] Suppose  $\eta\in \co_2$ satisfies 
$\ord_{k_v}(\tr_{k_2/k_v}(\eta)-a_1)=i$,   
$\n_{k_2/k_v}(\eta)\equiv a_2 \; (\gp_v^{i+1})$.  
If $\eta'\in \co_2$, 
$\tr_{k_2/k_v}(\eta')\equiv a_1 \; (\gp_v^i)$
and $\n_{k_2/k_v}(\eta')\equiv a_2 \; (\gp_v^{i+1})$, 
we have $\ord_{k_v}(\tr_{k_2/k_v}(\eta')-a_1)=i$. 
\end{itemize}
\end{prop}
\begin{proof} 
We first consider the case $i=\ell=1$.
For any uniformizer $\eta\in \co_2$,  
$\tr_{k_2/k_v}(\eta)\equiv a_1\equiv 0\;(\gp_v)$.  
So the statement (1) follows from the fact that 
any unit in $\co_v$ is a square modulo $\gp_v$.   
Consider (2).  By Lemma \ref{etaform}, there exist 
$e,f\in \co_v$ such that $\eta'=e\pi_v+(1+f\pi_v)\eta$. 
Let $\tr_{k_2/k_v}(\eta)=a_1+ h\pi_v$ with $h\in \co_v^{\times}$.  
Then 
\begin{align*}
\tr_{k_2/k_v}(\eta')
& = 2\pi_v e+(1+f\pi_v)\tr_{k_2/k_v}(\eta) \\
& = 2\pi_v e+(1+f\pi_v)(a_1+h\pi_v) \\
& \equiv a_1+h\pi_v  \; (\gp_v^2)
.\end{align*}
So $\ord_{k_v}(\tr_{k_2/k_v}(\eta')-a_1)=1$ also.  
This proves the proposition when $i=\ell=1$.  

Suppose $i\geq 2$ and  $\eta\in S_{i,i}(k_1,k_2)$.  
Let $\tr_{k_2/k_v}(\eta)-a_1=\gam_1 \pi_v^i,\; 
\n_{k_2/k_v}(\eta)-a_2= \gam_2\pi_v^i $, where 
$\gam_i\in \co_2$ for $i=1,2$.  
For (1), we look for an element of the form 
$\eta'=e\pi_v +(1+ f\pi_v)\eta$ with $e,f\in\co_v$.
If $\eta'$ satisfies the condition of (2), 
$\eta'$ is of the above form by Lemma \ref{etaform}.
Therefore, in both cases we consider $\eta'$ of the above
form.   Then 
\begin{equation*} 
\begin{aligned} 
\tr_{k_2/k_v}(\eta') 
& = 2e\pi_v +(1+ f\pi_v )(a_1+ \gam_1\pi_v^i), \\
\n_{k_2/k_v}(\eta') 
& = e^2\pi_v^2 +e\pi_v(1+ f\pi_v)(a_1+ \gam_1\pi_v^i) \\
& \quad + (1+f \pi_v)^2(a_2+\gam_2\pi_v^i)
.\end{aligned}
\end{equation*}
So 
\begin{equation}  \label{congform} 
\begin{aligned} 
\tr_{k_2/k_v}(\eta') -a_1
& \equiv 2e\pi_v +a_1f\pi_v +\gam_1\pi_v^i \; (\gp_v^{i+1}), \\
\n_{k_2/k_v}(\eta')  -a_2
& \equiv e^2\pi_v^2+a_1e\pi_v+a_1ef\pi_v^2 \\
& \quad +2a_2f\pi_v+a_2f^2\pi_v^2 +\gam_2\pi_v^i \; (\gp_v^{i+1})
.\end{aligned}
\end{equation}

Consider the case $2\leq i=\ell\leq m_v$.
We have 
\begin{equation}
\begin{aligned} 
\tr_{k_2/k_v}(\eta')-a_1 
& \equiv \gam_1\pi_v^i \; (\gp_v^{i+1}), \\
\n_{k_2/k_v}(\eta')-a_2 
& \equiv e^2\pi_v^2+a_2f^2\pi_v^2+\gam_2\pi_v^i \; (\gp_v^{i+1})
.\end{aligned}
\end{equation}
Since the orders of $e^2\pi_v^2$, $a_2f^2\pi_v^2$
can be any even or odd integer greater than or equal to 
two, we can choose $e,f$ so that
$\n_{k_2/k_v}(\eta') -a_2 \equiv 0 \;(\gp_v^{i+1})$. 
By the first congruence, we still have 
$\tr_{k_2/k_v}(\eta')-a_1 \equiv 0 \;(\gp_v^i)$. 
This proves (1).  If $\ord_{k_v}(\tr_{k_2/k_v}(\eta)-a_1)=i$, 
$\ord_{k_v}(\tr_{k_2/k_v}(\eta')-a_1)=i$ 
by the first congruence also. So this proves (2).  

Consider the case $2\leq i=\ell = m_v+1$.   
We have 
\begin{equation} 
\begin{aligned} 
\tr_{k_2/k_v}(\eta') -a_1
& \equiv 2e\pi_v +\gam_1\pi_v^i \; (\gp_v^{i+1}), \\
\n_{k_2/k_v}(\eta') -a_2
& \equiv e^2\pi_v^2 +a_2f^2\pi_v^2 +\gam_2\pi_v^i\;(\gp_v^{i+1})
.\end{aligned}
\end{equation}
As long as $e\in \co_v$, 
$\tr_{k_2/k_v}(\eta') -a_1\equiv 0 \; (\gp_v^i)$.  
By the same consideration as the previous case, 
we can choose $e,f$ so that
$\n_{k_2/k_v}(\eta') -a_2 \equiv 0 \;(\gp_v^{i+1})$.
This proves (1). 
We now turn to (2). By assumption, 
$\gam_2\pi_v^i\equiv  \n_{k_2/k_v}(\eta') -a_2 
\equiv 0 \; (\gp_v^{i+1})$.
So $e^2\pi_v^2 +a_2f^2\pi_v^2 \equiv 0 \; (\gp_v^{i+1})$.
Since the orders of $e^2\pi_v^2,a_2f^2\pi_v^2$ are even and odd, 
$e^2\pi_v^2,a_2f^2\pi_v^2\equiv 0 \; (\gp_v^{i+1})$.  
Since $i+1\geq 3$, 
$e\in \gp_v$.  So $2e\pi_v\equiv 0 \; (\gp_v^{i+1})$.  
This implies that $\ord_{k_v}(\tr_{k_2/k_v}(\eta')-a_1)=i$ which
proves (2).  

We now assume $\ell<i$.  Since $i<2\ell$ by assumption, 
$\ell>1$.  We first consider the case $\ell\leq m_v$.  
Then $\tr_{k_2/k_v}(\eta') -a_1 \equiv 0\; (\gp_v^i)$ 
if and only if there exists $h\in \co_v$ such that 
$f = -2e/a_1+ h\pi_v^{i-\ell-1}$.  Then 
by (\ref{congform}), 
\begin{equation} \label{congform2}
\begin{aligned} 
\tr_{k_2/k_v}(\eta') -a_1 
& \equiv (a_1/\pi_v^\ell) h \pi_v^i 
+ \gam_1 \pi_v^i \; (\gp_v^{i+1}), \\
\n_{k_2/k_v}(\eta') -a_2 
& \equiv e^2\pi_v^2+a_1e\pi_v+a_1e\pi_v^2(-2e/a_1+h\pi_v^{i-\ell-1})\\
& \quad +\gam_2\pi_v^i+2a_2\pi_v(-2e/a_1+h\pi_v^{i-\ell-1}) \\
& \quad +a_2\pi_v^2(-2e/a_1+h\pi_v^{i-\ell-1})^2 \\
& \equiv (-1+4a_2/a_1^2)(e^2\pi_v^2 - a_1e\pi_v) \\
& \quad + a_2h^2\pi_v^{2(i-\ell)} + 2a_2h \pi_v^{i-\ell}+ \gam_2\pi_v^i 
\;(\gp_v^{i+1})
.\end{aligned}
\end{equation}

Let $N_1= \ord_{k_v}(e),\; N_2 = \ord_{k_v}(h)$.  
Consider (1).  
We choose $e,h$ so that $0\leq N_1<\ell-1$ 
and $0\leq N_2< m_v-i+\ell$. This is possible because $\ell>1$.  
Then $\ord_{k_v}(e^2\pi_v^2) < \ord_{k_v}(a_1e\pi_v)$ and 
$\ord_{k_v}(h^2\pi_v^{2(i-\ell)}) < \ord_{k_v}(2h\pi_v^{i-\ell})$.  
Note that $\ord_{k_v}(e^2\pi_v^2) = 2N_1+2<2\ell$ and it can be any even 
integer between $2$ and $2\ell-2$.  Also 
$\ord_{k_v}(a_2h^2\pi_v^{2(i-\ell)}) < 2m_v+1$ and it can be any odd 
integer between $2(i-\ell)+1$ and $2m_v-1$.  
Since $i<2\ell$, $2(i-\ell)+1\leq i$.  
So we can choose $e,h$ so that 
$\n_{k_2/k_v}(\eta') -a_2 \equiv 0 \; (\gp_v^{i+1})$.  
Since $h\in\co_v$, the condition $\tr_{k_2/k_v}(\eta') -a_1 
\equiv 0 \; (\gp_v^i)$ is still satisfied.  This proves (1).  

Consider (2).  
If $N_2\geq m_v-i+\ell$ then $N_2>0$ because 
$m_v-i+\ell\geq 2\ell-i>0$, by assumption.  So $h\in \gp_v$.  
Therefore, $\ord_{k_v}(\tr_{k_2/k_v}(\eta') -a_1) = i$.  
So we assume that $N_2< m_v-i+\ell$.  If $N_1\geq \ell-1$
then $e^2\pi_v^2-a_1e\pi_v \in \gp_v^{2\ell}\sub \gp_v^{i+1}$
and so $a_2h^2\pi_v^{2(i-\ell)}\in \gp_v^{i+1}$.  
If $N_1<\ell-1$ then 
$e^2\pi_v^2 + a_2h^2\pi_v^{2(i-\ell)}\in \gp_v^{i+1}$. 
Since the orders of these elements are even and odd, 
$a_2h^2\pi_v^{2(i-\ell)}\in \gp_v^{i+1}$.  In both 
cases, $h^2\in \gp_v^{i-2(i-\ell)} = \gp_v^{2\ell-i}$.  
Since $2\ell-i>0$, 
$h\in \gp_v$.  Therefore, $\ord_{k_v}(\tr_{k_2/k_v}(\eta') -a_1) = i$. 
This proves (2). 

We now assume $\ell=m_v+1$ and so $i\leq 2m_v$.  
Then by (\ref{congform}), 
$\tr_{k_2/k_v}(\eta') -a_1 \equiv 0 \; (\gp_v^i)$ 
if and only if there exists $h\in \co_v$ such that 
$e= -(a_1/2)f+ h \pi_v^{i-m_v-1}$.  Then 
\begin{equation}  \label{congform4}
\begin{aligned}
\tr_{k_2/k_v}(\eta') -a_1 
& \equiv  (2/\pi_v^{m_v}) h\pi_v^i +\gam_1 \pi_v^i \;(\gp_v^{i+1}), \\
\n_{k_2/k_v}(\eta') -a_2
& \equiv ((-a_1/2)f+h\pi_v^{i-m_v-1})^2\pi_v^2 \\
& \quad + a_1\pi_v(1+f\pi_v)((-a_1/2)f+h\pi_v^{i-m_v-1})\\
& \quad + 2a_2f\pi_v+a_2f^2\pi_v^2+\gam_2\pi_v^i \\
& \equiv (a_2-(a_1^2/4))f^2\pi_v^2 + (2a_2-(a_1^2/2))f\pi_v \\
& \quad + h^2\pi_v^{2(i-m_v)} + \gam_2\pi_v^i
\;(\gp_v^{i+1})  
.\end{aligned}
\end{equation} 
Let $a_2-(a_1^2/4) = r\pi_v$ and $2a_2-(a_1^2/2)= s\pi_v^{m_v+1}$.  
Then it is easy to see that $r,s\in\co_v^{\times}$.  

Suppose $N=\ord_{k_v}(f)$.  
Consider (1).  We choose $0\leq N<m_v-1$.  Then 
$\ord_{k_v}(rf^2\pi_v^3) = 2N+3< N+m_v+2= \ord_{k_v}(sf\pi_v^{m_v+2})$ 
and  $\ord_{k_v}(rf^2\pi_v^3) = 2N+3$ can be any odd integer between 
$3$ and $2m_v-1$.  Since $i\leq 2m_v$, $2(i-m_v)-i = i-2m_v\leq 0$.  
So $\ord_{k_v}(h^2\pi_v^{2(i-m_v)})$ 
can be any even integer greater than or equal to $i$.  
Therefore, we can choose $f,h$ so that 
$\n_{k_2/k_v}(\eta') -a_2\equiv 0 \;(\gp_v^{i+1})$.  
This proves (1). 

Consider (2).  By assumption, 
\begin{equation*}
rf^2\pi_v^3+ sf\pi_v^{m_v+2} 
+ h^2\pi_v^{2(i-m_v)} \equiv 0 \;(\gp_v^{i+1})
.\end{equation*}
If $N\geq m_v-1$ then 
\begin{equation*}
rf^2\pi_v^3+ rf\pi_v^{m_v+2} \in \gp_v^{2m_v+1}\sub \gp_v^{i+1}
.\end{equation*}
So $h^2\pi_v^{2(i-m_v)} \equiv 0 \;(\gp_v^{i+1})$. 
If $N< m_v-1$ then 
$\ord_{k_v}(rf^2\pi_v^3)< \ord_{k_v}(rf\pi_v^{m_v+2})$ and 
the orders of $rf^2\pi_v^3,h^2\pi_v^{2(i-m_v)}$ 
are odd and even respectively. This implies that 
$h^2\pi_v^{2(i-m_v)} \equiv 0 \;(\gp_v^{i+1})$ also.  
In both cases, $h^2\in \gp_v^{2m_v+1-i}\sub \gp_v$ and so 
$h\in \gp_v$.   
Therefore,  $\ord_{k_v}(\tr_{k_2/k_v}(\eta') -a_1) = i$.  
This proves (2).  
\end{proof}

The following corollary is easily deduced from the proposition.
\begin{cor}\label{levelcond} 
The level of $k_1,k_2$ is $i$ if and only if 
there exists $\eta\in \co_2$ such that 
\begin{equation*}
\ord_{k_v}(\tr_{k_2/k_v}(\eta) - a_1)=i,\;
\n_{k_2/k_v}(\eta)\equiv a_2 \; (\gp_v^{i+1})
.\end{equation*}
\end{cor}

As we discussed in section 7 of \cite{kable-yukie-pbh-I},
the following 
proposition provides a relation between the level 
and the relative discriminant of $k_1\cdot k_2/k_2$.    
\begin{prop}\label{compare-discriminants}
\begin{itemize}
\item[(1)]  Suppose $\ell_1=\ell_2=\ell$, 
$\del_1=\del_2=\del$ and $\ell\leq \lev(k_1,k_2)\leq\del$. Let
$i=\lev(k_1,k_2)$. 
Then we have $\Del_{k_1\cdot k_2/k_2} = \gp_2^{2(\del-i)}$
and $\Del_{k_3/k_v} = \gp_v^{2(\del-i)}$.  
\item[(2)]  Suppose $\ell_1>\ell_2$. 
Then we have $\Del_{k_1\cdot k_2/k_1} 
= \gp_1^{\del_2},\; \Del_{k_1\cdot k_2/k_2}
= \gp_2^{2\del_1-\del_2}$. 
\end{itemize}
\end{prop}
\begin{proof} Consider the first claim in (1). If
$\lev(k_1,k_2)=\del$ then $k_1\cdot k_2/k_2$ and $k_3/k_v$ are
unramified, by Proposition \ref{unramified-char}, and so (1) holds. We
may now assume that $i<\del$. 
We choose $\eta$ which satisfies the condition of 
Corollary \ref{levelcond}.  Let $p(z)=z^2+a_1z+a_2=0$ be the 
Eisenstein equation with roots $\al=\{\al_1,\al_2\}$, 
which generate the field $k_1$.  Let 
$\gam=\pi_2^{-i}(\al_1+\eta)$.  Then 
$p(\al_1)=0$ is equivalent to the following equation
\begin{equation}  \label{neweisen}
\gam^2 + \pi_2^{-i}(a_1-2\eta)\gam 
+ \pi_2^{-2i}(\eta^2-a_1\eta+a_2) =0
.\end{equation}

Since $\eta^2-a_1\eta+a_2=\eta(\tr_{k_2/k_v}(\eta)-a_1)
+a_2-\n_{k_2/k_v}(\eta)$, the order of the 
third term in (\ref{neweisen}) is one.  Therefore, 
(\ref{neweisen}) is an Eisenstein equation whose roots
generate $k_1\cdot k_2/k_2$.  
If $\ell\leq m_v$ then $\ord_{k_2}(a_1)=2\ell$ and 
$\ord_{k_2}(2\eta)=2m_v+1>2\ell$.  
So $\ord_{k_2}(\pi_2^{-i}(a_1-2\eta)) = \del-i$.  
If $\ell=m_v+1$ then $\ord_{k_2}(a_1)=2m_v+2$ and 
$\ord_{k_2}(2\eta)=2m_v+1=\del$.  
So $\ord_{k_2}(\pi_2^{-i}(a_1-2\eta)) = \del-i$ also.   
Since $2\co_2=\gp_2^{2m_v}$ and 
$\del-i\leq 2m_v$, $\Del_{k_1\cdot k_2/k_2} = \gp_2^{2(\del-i)}$. 

Consider the second claim in (1).  
By Corollary \ref{levelcond}, we choose an 
element $\eta\in S_{i,i+1}(k_1,k_2)$.  
We may assume that $-\eta$ is one of the  roots of $p_2(z)$.  
Let $p_3(x)$ be the polynomial (\ref{p3defn}).  
Then the roots of $p_3(z)$ generate the field $k_3$.  
We evaluate the order of the element $J$ in 
Lemma \ref{k3generator}, which is the same as that in 
(\ref{J-value}).    

By assumption $\ord_{k_v}(a_2-b_2)\geq i+1$, 
$\ord_{k_v}(a_1-b_1)=i$, 
$\ord_{k_v}(a_1(b_2-a_2))\geq i+2$ and 
$\ord_{k_v}(a_2(a_1-b_1))=i+1$.  
Therefore, $\ord_{k_v}(J) = 2i+1$. 
Now 
\begin{equation}  \label{eisen-k3}
\pi_v^{-2i}p_3(\pi_v^i z)
= z^2 - \pi_v^{-i}[2(a_2+b_2)-a_1b_1]z+\pi_v^{-2i}J
.\end{equation}
Note that $2(a_2+b_2)=4a_2+2(b_2-a_2)$
and $\ord_{k_v}(4a_2) = 2m_v+1,\; 
\ord_{k_v}(2(b_2-a_2)) \geq m_v+i+1$.  
If $\ell\leq m_v$ then $\ord_{k_v}(a_1b_1) = 2\ell
< 2m_v+1,m_v+i+1$ since $i\geq \ell$.  
This implies that the order of 
$\pi_v^{-i}[2(a_2+b_2)-a_1b_1]$ is $2\ell-i=\del-i$.  
If $\ell=m_v+1$ then $\ord_{k_v}(a_1b_1)\geq 2m_v+2=\del+1$.  
Since $i\geq m_v+1$, the order of 
$2\pi_v^{-i}(a_2+b_2)$ is $2m_v+1-i=\del-i$.  
Therefore, in both cases, (\ref{eisen-k3}) 
is an Eisenstein polynomial with the order 
of the coefficient of the 
middle term $\del-i\leq m_v$.   Therefore, 
$\Del_{k_3/k_v} = \gp_v^{2(\del-i)}$.   

Consider (2).  By Lemma \ref{bqbound}, 
$\del_3=\del_1>\del_2$.  Let 
$i = \lev(k_1,k_3)$.  Then 
$\del_2 = 2(\del_1-i)$ by the second 
statement of  (1).  Therefore, using the 
first statement of (1), $\Del_{k_1\cdot k_2/k_v}
= \Del_{k_1\cdot k_3/k_v} = \gp_v^{2\del_1+\del_2}$
(see \cite{weilc}, Corollary 4, p. 142 which is a 
local version of \cite{weilc}, Proposition 13, p. 156).  
This implies that $\Del_{k_1\cdot k_2/k_1} 
=\gp_1^{\del_2}$ and $\Del_{k_1\cdot k_2/k_2}
=\gp_2^{2\del_1-\del_2}$.  
Thus (2).  
\end{proof}  

We now review the equivalence relation $x\asymp y$
and explain the notation in the introduction.   
Since we are only concerned with $x$ such that $k_v(x)/k_v$ is ramified,
we restrict ourselves to such orbits.  
Suppose $x,y\in V^{\sst}_{k_v}$ and $k_v(x),k_v(y)$
are ramified  quadratic extensions of $k_v$.  
If the type of $x$ is (rm~rm)*
or (rm~rm~ur), $x\asymp y$ means $x,y$ are in the
same $G_{k_v}$-orbit.  If the type of $x$ is 
(rm~rm~rm), we write $x\asymp y$ if and only if 
$\Del_{k_v(x)/k_v}=\Del_{k_v(y)/k_v}$ and 
$\lev(k_v(x),\ti k_v) = \lev(k_v(y),\ti k_v)$.  
By Proposition \ref{compare-discriminants}, 
the last condition is equivalent to 
the condition $\Del_{\ti k_v(x)/\ti k_v}
=\Del_{\ti k_v(y)/\ti k_v}$.  
If $k_v(x)/k_v$ is ramified, we let 
$\Del_{k_v(x)/k_v} =\gp_v^{\del_{x,v}}$
and $\lam_{x,v}=\lev(k_v(x),\ti k_v)$.  
This explains the notation in Tables \ref{table-dyadic-ungrouped}
and \ref{table-dyadic-grouped}.  

As we promised earlier, we explain our motivation for 
our formulation.  
Before finally choosing the formulation of the filtering
process in section 6 of \cite{kable-yukie-pbh-I},
we carried out some experiments.  
At first we tried to compute the standard local zeta functions
explicitly and we did succeed for non-dyadic places, even though
we later settled on a uniform estimate without the explicit forms
to shorten the paper.  
Then we worked on dyadic places and we discovered 
that it is difficult even to determine the constant terms
of the standard local zeta functions.  If one tries to 
compute them, 
the set $S_{i_1,i_2}(k_1,k_2)$ naturally occurs.
In fact, if $v$ is dyadic and an orbit of $x\in V^{\sst}_{k_v}$ 
corresponds  to a field $k_v(x)$ such that 
$\lev(k_v(x),\ti k_v)=i$, it turns out that the constant term 
of the standard local zeta function is 
$\sum_{j=0}^i \big(\nn_1(k_v(x),\ti k_v,j)
+ \nn_2(k_v(x),\ti k_v,j)\big)$.  We also evaluated the terms in this
sum and the answer was that if $x$ is of type (rm~rm~rm) then
$\nn_r(k_v(x),\ti k_v,j)= q_v^j$ for 
$r=1,2$ and $j\leq i$ and if $x$ is of type (rm~rm~ur)
then $\nn_r(k_v(x),\ti k_v,j)= q_v^j$ for $r=1,2$ 
and $j\leq \ti\del-1$ and 
$\nn_1(k_v(x),\ti k_v,\ti\del)= q_v^{\ti\del},\; 
\nn_2(k_v(x),\ti k_v,\ti\del)=0$.
In the process we had to prove something like 
Proposition \ref{crucial-prop}.  We realized later 
that we did not need the 
constant term  nor any estimate of the 
standard local zeta functions at dyadic places, but 
having Proposition \ref{crucial-prop} eventually helped
us to evaluate the local densities at dyadic places.  
This was our motivation for introducing the set 
$S_{i_1,i_2}(k_1,k_2)$.  

\section{The volume of the integral points of the stabilizer}
\label{stab-volume}

In this section we evaluate $\vol(K_v\cap G_{x\,k_v}^{\circ})$
for orbits of types (rm~rm~ur) and (rm~rm~rm).  
The measure on $G_{x\,k_v}^{\circ}$ is defined in Definition 5.13
of \cite{kable-yukie-pbh-I}. 
We shall not repeat the definition here but 
instead recall an alternative formula.  
Consider the usual multiplicative measure on $\ti k_v(x)$, i.e.
that for which the volume of $\co_{\ti k_v(x)}^{\times}$ is $1$.  
Suppose $\al_1$ is a uniformizer of $k_v(x)$. 
Then it  was proved in Lemma 10.4 of \cite{kable-yukie-pbh-I} that 
\begin{equation} \label{alternative}
\vol(K_v\cap G_{x\,k_v}^{\circ})
= \vol(\ti\co_v[\al_1]^{\times})
.\end{equation}

To determine the volume in the above cases 
we need the following result,
which provides a reinterpretation of the level in these cases.

\begin{prop}
Suppose that $k_1,k_2/k_v$ are distinct ramified quadratic
extensions. Let $\co_i$ be the integer ring of $k_i$, $\gp_i$ be the
prime ideal of $\calo_i$ and 
$\pi_i$ be its uniformizer for $i=1,2$.  We denote 
the integer ring of $k_1\cdot k_2$ by $\co_{k_1\cdot k_2}$.  
Let $p_1(z)=z^2+a_1z+a_2$ be an Eisenstein polynomial
defining $k_1$ and $\al_1\in\co_1$ be a root of $p_1$. Let $f$ be
the least integer such that $\gp_2^f\cdot
\co_{k_1\cdot k_2}\subseteq\co_2[\al_1]$. 
Then $f=\lev(k_1,k_2)$.
\end{prop}
\begin{proof}
We shall show first that $f\geq\lev(k_1,k_2)$ in general.
Let $i=\lev(k_1,k_2)$ and choose $\beta_1\in\co_2$ with minimal
polynomial $z^2+b_1z+b_2$ such that $a_1\equiv b_1\; (\gp_v^i)$ and
$a_2\equiv b_2\; (\gp_v^i)$. Consider the element
$(\al_1-\beta_1)/\pi_2^i$ of $k_1\cdot k_2$; we claim that it is an
integer. In fact, by Lemma \ref{k3generator}, 
\begin{align*}
&\operatorname{N}_{k_1\cdot k_2/k_v}((\al_1-\beta_1)/\pi_2^i) \\
&=\frac{1}{\n_{k_2/k_v}(\pi_2)^{2i}}
[(a_2-b_2)^2+(a_1-b_1)(a_1b_2-a_2b_1)] \\
&=\frac{1}{\n_{k_2/k_v}(\pi_2)^{2i}}
[(a_2-b_2)^2+(a_1-b_1)(a_1(b_2-a_2)+
a_2(a_1-b_1))]\in\co_v\,,
\end{align*}
by hypothesis. Thus $\pi_2^{f-i}(\al_1-\beta_1)\in\co_2[\al_1]$ 
and it follows that $f\geq i$, as claimed. 

We know that $\lev(k_1,k_2)\geq1$ and so $f\geq 1$ and
$\lev(k_1,k_2)=f$ if $f=1$. We now assume that $f\geq2$ to complete
the proof.
There are $\eta,\zeta\in\co_2$ such
that $(\eta+\zeta\al_1)/\pi_2^f\in\co_{k_1\cdot k_2}$ and one of
$\eta,\zeta$ is a unit (for otherwise $f$ would not be the least
integer with its defining property). Taking norms from $k_1\cdot
k_2$ to $k_2$ we find that
$\eta^2-a_1\eta\zeta+a_2\zeta^2\in\gp_2^{2f}$. 
Since $a_1,a_2\in \gp_v\sub \gp_2^2$, $\eta\in\gp_2$.
It follows that $\zeta\in\co_2^{\times}$. Furthermore,
$\eta^2+a_2\zeta^2\in\gp_2^3$ and
$a_2\zeta^2\in\gp_2^2\setminus\gp_2^3$ from which it follows that
$\eta\in\gp_2\setminus\gp_2^2$. Let us set $\varpi=\eta/\zeta$.
Then $\varpi$ is a uniformizer of $k_2$ and
$\varpi^2-a_1\varpi+a_2\in\gp_2^{2f}$. Let $z^2-c_1z+c_2$ be the
(Eisenstein) minimal polynomial of $\varpi$ over $k_v$. Then
\begin{align*}
(c_1-a_1)\varpi+(a_2-c_2)&=(\varpi^2-a_1\varpi+a_2)-
(\varpi^2-c_1\varpi+c_2) \\
&=(\varpi^2-a_1\varpi+a_2)\in\gp_2^{2f}\,.
\end{align*}
Since $(c_1-a_1)\varpi$ has odd order in $k_2$ and $(a_2-c_2)$ has
even order, it follows that $(a_2-c_2)\in\gp_2^{2f}
\cap\co_v=\gp_v^{f}$ and $(c_1-a_1)\varpi\in\gp_2^{2f+1}$ which
implies that $(c_1-a_1)\in\gp_v^{f}$. Thus $f\leq\lev(k_1,k_2)$.
This proves the proposition.  
\end{proof}
\begin{prop}  \label{stabvolume}
If $x$ is the standard orbital representative for an orbit with
type {\upshape (rm rm ur)} then $\vol(K_v\cap G_{x\,k_v}^{\circ})=
(1+q_v^{-1})^{-1}q_v^{-\ti\del_v}$.   
If $x$ is the standard orbital representative for an orbit with
type {\upshape (rm rm rm)} then $\vol(K_v\cap G_{x\,k_v}^{\circ})=
q_v^{-i}$ where $i=\lev(k_v(x),\kt_v)$.
\end{prop}
\begin{proof}
The ring $\ti\co_v[\al_1]$ is an $\ti\co_v$-order in
$\co_{\kt_v(x)}$ and so if $\beta_1\in\co_{\kt_v(x)}$ satisfies
$\co_{\kt_v(x)}=\ti\co_v[\beta_1]$ then there is some $i\geq0$ such
that
\begin{equation*}
\ti\co_v[\al_1]=\{a+b\beta_1\mid a\in\ti\co_v,b\in\ti\gp_v^{i}\}\,.
\end{equation*}
From the previous proposition we see that
$i=\lev(k_v(x),\kt_v)\geq1$. Then
\begin{equation*}
\ti\co_v[\al_1]^{\times}=\{a+b\beta_1\mid a\in\ti\co_v^{\times},
b\in\ti\gp_v^i\}\,.
\end{equation*}
The normalized additive Haar measure on $\co_{\kt_v(x)}$ is
$da\,db$ and so the normalized multiplicative Haar measure on
$\co_{\kt_v(x)}^{\times}$ is $(1-q_{\ti k_v(x)}^{-1})^{-1}da\,db$,
where $q_{\ti k_v(x)}$ is the module of $\kt_v(x)$. 
Since $\ti\co_v/\ti\gp_v\cong \co_v/\gp_v$, 
\begin{equation*}
\vol(\ti\co_v[\al_1]^{\times})
=(1-q_{\ti k_v(x)}^{-1})^{-1}(1-q_v^{-1})q_v^{-i}\,.
\end{equation*}
In case the index is {\upshape (rm rm ur)}, 
$q_{\ti k_v(x)}=q_v^2$ and $i=\ti\del_v$ and we have
$\vol(\ti\co_v[\al_1]^{\times})=(1+q_v^{-1})^{-1}q_v^{-\ti\del_v}$. 
In case the index is {\upshape (rm rm rm)}, $q_{\ti k_v(x)}=q_v$ and
we have $\vol(\ti\co_v[\al_1]^{\times})=q_v^{-i}$. 
\end{proof}

\section{Orbital volumes at the ramified dyadic places} 
\label{bxrf} 
 
In this section, we group orbits according to the 
level and compute  $\sum_x\vol(K_v x)$ for each group of orbits. 
 
Let $\ti p(z)= z^2+b_1z+b_2$ be an Eisenstein polynomial whose
roots $\eta=\{\eta_1,\eta_2\}$  generate $\ti k_v$.  
Let $\ti\ell=\ord_{k_v}(b_1)$ if 
$\ord_{k_v}(b_1)\leq m_v$ and $\ti\ell=m_v+1$ if 
$\ord_{k_v}(b_1)\geq m_v+1$, and 
$\ti\del_v =2\ti\ell$ or $2m_v+1$, as before.  
For an Eisenstein polynomial  $p_1(z)=z^2+a_1z+a_2 =0$
we define $\ell(p_1)$ and $\del(p_1)$ similarly.  
\begin{defn} \label{xtype}
\begin{itemize}
\item[(1)] If $\ell_1 \not=\ti\ell$ then $\gX_{\ell_1}$ is the set of 
isomorphism classes of quadratic extensions $k'$ of $k_v$ generated 
by roots of an Eisenstein equation $p_1(z)=z^2+a_1z+a_2 =0$
such that $\ell(p_1)=\ell_1$.  
\item[(2)] If $\ti\ell\leq i<\ti\del_v$ then
$\gX_{\ti\ell}(i)$ is the set of isomorphism classes of 
quadratic extensions $k'$ of $k_v$ generated 
by roots of an Eisenstein equation $p_1(z)=z^2+a_1z+a_2 =0$
such that $\ell(p_1)=\ti\ell$ and $\lev(k',\ti k_v)=i$.  
\item[(3)] We define $\gX^{\text{\upshape ur}}_{\ti\ell}$ 
to be the singleton containing the unique quadratic extension 
of $k_v$ of type (rm rm ur).  
\item[(4)] We define $\gX^*_{\ti\ell}$ to 
be the singleton containing the unique quadratic extension 
of $k_v$ of type (rm rm)*. 
\end{itemize}
\end{defn}  

Let $a_0(x),a_1(x),a_2(x)$ be as in 
(\ref{aform}).  
For each type of $x$ in Definition \ref{xtype}, we compute 
$\sum_x \vol(K_v x)$.  Our strategy is the same as that in 
section 11 of \cite{kable-yukie-pbh-I};
we define a subset ${\mathcal D}\sub V_{\co_v}$
using congruence conditions, cover $K_vx$ by disjoint copies
of ${\mathcal D}$ and count the number of copies.  Our first task 
is to define the set ${\mathcal D}$ for each case, which we shall do
as follows.  

Let $\ell_1\not=\ti\ell$.  We put $\ell=\operatorname{min}
\{\ell_1,\ti\ell\}$.  We define
${\mathcal D}_{\ell_1}$ to be the set of $x$ which satisfy
the conditions
\begin{equation}\label{lcond}
\begin{aligned}
{} & x_{11}\in \ti\co_v^{\times},\; x_{20}\in\co_v^{\times}, \\
&  \ord_{\ti k_v}(x_{21})=1,\; \ord_{k_v}(x_{12})=\ell,\;  
\ord_{k_v}(x_{22})\geq\ell+1, \\
&  \ord_{k_v}(a_1(x))
\begin{cases} = \ell_1 &\text{\quad if\ } \ell_1\leq m_v, \\ 
\geq m_v+1 &\text{\quad if
\ } \ell_1=m_v+1.\end{cases}
\end{aligned}
\end{equation} 
We define ${\mathcal D}_{\ti\ell}(i)$ 
to be the set of $x$ which satisfy
the condition
\begin{equation}\label{licond}
\begin{aligned}
{} & x_{11}\in \ti\co_v^{\times},\; x_{20}\in\co_v^{\times}, \\
& \ord_{\ti k_v}(x_{21})=1,\; \ord_{k_v}(x_{12})=i,\; 
\ord_{k_v}(x_{22})\geq i+1, \\
&  \ord_{k_v}(a_1(x))
\begin{cases} = \ti\ell &\text{\quad if\ } \ti\ell\leq m_v, \\ 
\geq m_v+1 &\text{\quad if\ } \ti\ell=m_v+1.\end{cases}
\end{aligned}
\end{equation} 
We define  ${\mathcal D}^{\text{ur}*}_{\ti\ell}$ 
to be the set of $x$ which satisfy
the conditions
\begin{equation}\label{ur*cond}
\begin{aligned}
{} & x_{11}\in \ti\co_v^{\times},\; x_{20}\in\co_v^{\times}, \\
& \ord_{\ti k_v}(x_{21})=1,\; \ord_{k_v}(x_{12}),\, 
\ord_{k_v}(x_{22})\geq\ti\del_v, \\
&  \ord_{k_v}(a_1(x))
\begin{cases} = \ti\ell &\text{\quad if\ } \ti\ell\leq m_v, \\ 
\geq m_v+1 &\text{\quad if\ } \ti\ell=m_v+1.\end{cases}
\end{aligned}
\end{equation} 

Let $\eta=(\eta_1,\eta_2), \ti p(z)$ be as in the beginning 
of this section.  We define 
\begin{equation}\label{weta}
w_{\eta} = \left( \pmatrix 0 & 1\\ 1 & 0\endpmatrix,
\pmatrix 1 & -\eta_2\\ -\eta_1 & 0\endpmatrix \right)
.\end{equation}
Then $k_v(w_{\eta})= \ti k_v$.  Note that 
$w_{\eta} = (n(\eta_2),1)w_{\ti p}$.  
We define 
\begin{equation}\label{*cond}
{\mathcal D}^*_{\ti\ell} = \{x\in V_{\co_v}\mid x\equiv w_{\eta} 
\;(\gp_v^{\ti\del_v+1},\ti\gp_v^{2(\ti\del_v+1)})\}
.\end{equation}

Our next task is to show that points in the above sets correspond
to fields of types (1)--(4) in Definition \ref{xtype}.  
Given points in the above sets, we try to simplify them by group elements
as much as possible so that, after the simplification, the types of 
the corresponding fields are easy to determine.  
For this purpose we define
subgroups of $K_v$ which stabilize the above sets. 
They will also be used later for the computation of
$\sum_x\vol(K_vx)$. 
We use the coordinate system (\ref{gform}).   
\begin{defn} \label{hjdefn} 
For $j\geq 0$ we define 
\begin{align*}
H(j) & = \{g=(g_1,g_2)\in G_{\co_v}\mid g_{121}\in \ti\gp_v^{j+1},\;
g_{221}\in \gp_v\}, \\
\overline{H}(j) & = \{g=(g_1,g_2)\in 
\gl(2)_{\ti\co_v/\ti\gp_v^{j+1}} \times \gl(2)_{\co_v/\gp_v} 
\mid g_{121}=0,\; g_{221} = 0\}, \\
G(\pi_v^{\ti\del_v+1}) & = \{g=(g_1,g_2)\in G_{\co_v} \mid 
g_1\equiv 1\;(\ti\gp_v^{2(\ti\del_v+1)}),\; 
g_2\equiv 1\; (\gp_v^{\ti\del_v+1})\}
.\end{align*}
\end{defn} 

We put 
\begin{equation}\label{hjorder}
\begin{aligned}
Q(j) & = \# (K_v/H(j)) = \# (\gl(2)_{\ti\co_v/\ti\gp_v^{j+1}} 
\times \gl(2)_{\co_v/\gp_v})/\#\overline{H}(j) \\
& = \frac {(q_v^2-q_v)^2(q_v^2-1)^2q_v^{4j}}
{(q_v-1)^2q_v^{2j}q_v^{j+1}(q_v-1)^2q_v} \\
& = q_v^{j+2}(1+q_v^{-1})^2
.\end{aligned}
\end{equation}
\begin{prop}\label{hinvariant}
\begin{itemize}
\item[(1)] 
If $\ell_1 \not=\ti\ell$ and 
$\ell = \operatorname{min}\{\ell_1,\ti\ell\}$
then $H(\ell){\mathcal D}_{\ell_1} = {\mathcal D}_{\ell_1}$. 
\item[(2)] If $\ti\ell\leq i<\ti\del_v$ then 
$H(i){\mathcal D}_{\ti\ell}(i)={\mathcal D}_{\ti\ell}(i)$. 
\item[(3)] We have $H(\ti\del_v-1){\mathcal
D}_{\ti\ell}^{\text{\upshape ur}*}
={\mathcal D}_{\ti\ell}^{\text{\upshape ur}*}$.  
\item[(4)] We have $G(\pi_v^{\ti\del_v+1}){\mathcal D}_{\ti\ell}^{*}
= {\mathcal D}_{\ti\ell}^{*}$.  
\end{itemize}
\end{prop} 
\begin{proof} Part (4) is obvious.  Consider parts (1)--(3).
We put  $j=\ell,i,\ti\del_v-1$ for (1)--(3), respectively.
Then the group in question is $H(j)$ in all parts.

If $g=(g_1,g_2)\in H(j)$ then $(g_1,1),(1,g_2)\in H(j)$ and
$g=(1,g_2)(g_1,1)$. Thus it is enough to verify the claims for
$g=(g_1,1)$ and $g=(1,g_2)$ separately. We begin with $g=(g_1,1)$.
For $x$ in the form (\ref{xform}), let $y=gx=(y_1,y_2)$ and
consider similar coordinates for $y$. Then  
\begin{equation}\label{yform} 
\begin{aligned}
y_{r0} & = \n_{\ti k_v/k_v}(g_{111})x_{r0} 
+ \tr_{\ti k_v/k_v}(g_{111}g_{112}^{\sig}x_{r1}) 
+ \n_{\ti k_v/k_v}(g_{112})x_{r2}, \\
y_{r1} & = g_{111}g_{121}^{\sig}x_{r0}
+ g_{111}g_{122}^{\sig}x_{r1}
+ g_{112}g_{121}^{\sig} x_{r1}^{\sig} 
+ g_{112}g_{122}^{\sig}x_{r2}, \\
y_{r2} & = \n_{\ti k_v/k_v}(g_{121})x_{r0} 
+ \tr_{\ti k_v/k_v}(g_{121}g_{122}^{\sig}x_{r1}) 
+ \n_{\ti k_v/k_v}(g_{122})x_{r2}
\end{aligned}
\end{equation} 
for $r=1,2$. 

Suppose $x\in {\mathcal D}_{\ell_1}, {\mathcal D}_{\ti\ell}(i)$
or ${\mathcal D}_{\ti\ell}^{\text{ur}*}$.  Note that $j+1\geq 2$ 
in all cases.  So 
$g_{121},g_{121}^{\sig}\in\ti\gp_v^{j+1}\sub \ti\gp_v^2$,
$x_{12}\in\gp_v$ and $x_{11}\in \ti\co_v^{\times}$.  Therefore, 
$y_{11}\in\ti\co_v^{\times}$ by (\ref{yform}).  
We also have $y_{20}\in\co_v^{\times}$
since 
$\n_{\ti k_v/k_v}(g_{111})x_{20}\in\co_v^{\times},\; 
x_{21}\in\ti\gp_v$ and $x_{22}\in\gp_v$.
Since $g_{121},g_{121}^{\sig}\in\ti\gp_v^2,\; 
x_{22}\in\gp_v\sub \ti\gp_v^2$ and $\ord_{\ti k_v}(x_{21})=1$,
we further have $\ord_{\ti k_v}(y_{21})=1$.  

Note that $j+1\leq \ti\del_v$ in all cases.  
By Lemma \ref{jorder}, 
\begin{equation*}
\ord_{k_v}(\tr_{\ti k_v/k_v}(g_{121}g_{122}^{\sig}x_{r1}))
\geq \lfloor(j+1+\ti\del_v)/2\rfloor\geq j+1\,.
\end{equation*}
Therefore, 
$\tr_{\ti k_v/k_v}(g_{121}g_{122}^{\sig}x_{r1})\in \gp_v^{j+1}$
for $r=1,2$.  Also $\n_{\ti k_v/k_v}(g_{121})\in \gp_v^{j+1}$.
By assumption $\ord_{k_v}(\n_{\ti k_v/k_v}(g_{122})x_{12}) = j$
in cases (1), (2) and so $\ord_{k_v}(y_{12}) = j$.  In case (3), 
$\n_{\ti k_v/k_v}(g_{122})x_{12}\in \gp_v^{j+1}$
and so $y_{12}\in \gp_v^{j+1}$.  In all cases 
$\n_{\ti k_v/k_v}(g_{122})x_{22}\in \gp_v^{j+1}$
and so $y_{22}\in \gp_v^{j+1}$. 
We have $a_1(y)=\n_{\ti k_v/k_v}(\det(g_1))a_1(x)$ and
$\det(g_1)\in\ti\calo_v^{\times}$. Thus
$\ord_{k_v}(a_1(y))=\ord_{k_v}(a_1(x))$. 
All the conditions for $y$ to lie in $\mathcal D_{\ell_1}$, $\mathcal
D_{\ti\ell}(i)$ or $\mathcal D_{\ti\ell}^{\text{ur}*}$ have now been
verified and so 
$gx\in {\mathcal D}_{\ell_1},
{\mathcal D}_{\ti\ell}(i)$ or 
${\mathcal D}_{\ti\ell}^{\text{ur}*}$.  

We now assume that $g=(1,g_2)$. It is easy to verify that $g$
preserves all the conditions in (\ref{lcond})--(\ref{ur*cond}), with
the possible exception of the last. The necessary calculation to show
that this condition is also preserved by $g$ has already been carried
out in \cite{kable-yukie-pbh-I}, Lemma 11.23 and will not be repeated
here (note that $g_3$ stands in for $g_2$ in the proof of Lemma
11.23). This proves the proposition.
\end{proof} 
\begin{prop}\label{xbelong}
\begin{itemize}
\item[(1)] If $\ell_1\not=\ti\ell$ and
$x\in {\mathcal D}_{\ell_1}$ then
$k_v(x)\in \gX_{\ell_1}$.
If $k'\in \gX_{\ell_1}$ then there is some
$x\in {\mathcal D}_{\ell_1}$ with $k_v(x)\cong k'$ over $k_v$.
\item[(2)] If $\ti\ell\leq i<\ti\del_v$ and
$x\in {\mathcal D}_{\ti\ell}(i)$ then
$k_v(x)\in \gX_{\ti\ell}(i)$.
If $k'\in \gX_{\ti\ell}(i)$ then there is some
$x\in {\mathcal D}_{\ti\ell}(i)$ with $k_v(x)\cong k'$ over $k_v$.
\item[(3)] 
If $x\in {\mathcal D}^{\text{\upshape ur}*}_{\ti\ell}$ then
$k_v(x)\in \gX^{\text{\upshape ur}}_{\ti\ell}\cup \gX^*_{\ti\ell}$. If
$k'\in \gX^{\text{\upshape ur}}_{\ti\ell}\cup 
\gX^*_{\ti\ell}$ then there is some
$x\in {\mathcal D}^{\text{\upshape ur}*}_{\ti\ell}$ such that 
$k_v(x)\cong k'$ over $k_v$.
\item[(4)]  If $x\in {\mathcal D}^*_{\ti\ell}$ then $k_v(x)=\ti k_v$.  
\end{itemize}
\end{prop} 
\begin{proof}   
We first consider the first implication in each of (1)--(4).  
Let $x\in  {\mathcal D}_{\ell_1}, 
{\mathcal D}_{\ti\ell}(i)$, ${\mathcal D}^{ur*}_{\ti\ell}$,
or ${\mathcal D}^*_{\ti\ell}$.
Applying the element $g=(1,{}^tn(-x_{10}x_{20}^{-1}))$ to $x$, which
is permissible by
Proposition \ref{hinvariant}, we may assume that $x_{10}=0$; note that
this doesn't change $k_v(x)$. Further, applying
$g=(a(x_{11}^{-1},1),\n_{\ti k_v/k_v}(x_{11})a(1,x_{20}^{-1}))$ 
we may also assume that
$x_{11}=x_{20}=1$. This implies that $a_0(x)=1$ and
\begin{equation}\label{a1a2}
a_1(x) = \tr_{\ti k_v/k_v}(x_{21})-x_{12},\; 
a_2(x) =  \n_{\ti k_v/k_v}(x_{21})-x_{22}
.\end{equation} 
In case (4), $x_{22}\equiv 0 \;(\gp_v^{\ti\del_v+1})$ and
so $x_{22}\equiv 0 \;(\gp_v^2)$.
Note that $\ti\ell+1,i+1,\ti\del_v \geq 2$ for (1)--(3), respectively,  
and so $x_{22}\equiv 0 \;(\gp_v^2)$ in these cases also. Therefore,  
\begin{equation*} 
\ord_{k_v}(a_2(x))=\ord_{k_v}(\n_{\ti k_v/k_v}(x_{21})-x_{22})
= \ord_{k_v}(\n_{\ti k_v/k_v}(x_{21})) = 1
.\end{equation*}
By this and the last conditions in (\ref{lcond})--(\ref{ur*cond}), 
$F_x(z,1)$ is an  Eisenstein polynomial such that the 
corresponding $\ell$ is  
$\ell_1,\ti\ell,\ti\ell$ for (1)--(3),
respectively. 
In case (4), $x_{12}\equiv 0 \;(\gp_v^{\ti\del_v+1})$ and
so $\ord_{k_v}(a_1(x)) = \ti\ell$ or $\ord_{k_v}(a_1(x))\geq m_v+1$,
by Lemma \ref{tracelemma}.  
The first implication in (1) is now clear. 

Consider (2).  By assumption, 
\begin{equation*}  
\begin{aligned}
\ord_{k_v}(\tr_{\ti k_v/k_v}(x_{21}) -a_1(x)) 
& = \ord_{k_v}(x_{12})=i, \\ 
\n_{\ti k_v/k_v}(x_{21})-a_2(x) & = x_{22} \equiv 0 \; (\gp_v^{i+1})
.\end{aligned}
\end{equation*}
So, by Corollary \ref{levelcond}, $\lev(k_v(x),\ti k_v)=i$.  
This proves the first implication of (2).  

Consider (3). We have 
\begin{equation*}  
\begin{aligned}
\tr_{\ti k_v/k_v}(x_{21}) -a_1(x) & = x_{12}\equiv 0 \;
(\gp_v^{\ti\del_v}), \\ 
\n_{\ti k_v/k_v}(x_{21})-a_2(x) & = x_{22} \equiv 0 \; (\gp_v^{\ti\del_v})
.\end{aligned}
\end{equation*}
Therefore $S_{\ti\del_v,\ti\del_v}(k_v(x),\ti k_v)\not=\emptyset$. 
So the only possible types are (rm~rm~ur), (rm~rm)*,
by Corollary \ref{ib} and Proposition \ref{unramified-char}.  

By similar considerations, 
$S_{\ti\del_v+1,\ti\del_v+1}(k_v(x),\ti k_v)\not=\emptyset$
in case (4).  
So the only possible type is (rm rm)*, by the remark after
Corollary \ref{ib}. 

We now consider the second implication of (1)--(3). 
Suppose the roots of an Eisenstein equation $p(z)= z^2+a_1z+a_2=0$ 
generate $k'$. In part (1),
there exists $\eta\in \ti k_v$ such that 
$\ord_{k_v}(\tr_{\ti k_v}(\eta)-a_1)=\ell$ and 
$\n_{\ti k_v}(\eta)-a_2\equiv 0\; (\gp_v^{\ell+1})$,
by Proposition \ref{vanish}. In part (2),
by Corollary \ref{levelcond}, 
there exists $\eta\in \ti k_v$ such that 
$\ord_{k_v}(\tr_{\ti k_v}(\eta)-a_1)=i$ and 
$\n_{\ti k_v}(\eta)-a_2\equiv 0\; (\gp_v^{i+1})$. 
In part (3), by Proposition \ref{unramified-char}, 
there exists $\eta\in\ti k_v$ such that 
$\tr_{\ti k_v}(\eta)-a_1,\n_{\ti k_v}(\eta)-a_2\in\gp_v^{\ti\del_v}$. 
Let $x= (n(\eta-a_1),1)w_p$ in all cases. Then 
\begin{equation*}
x = \left( \pmatrix 0 & 1\\ 1 & \tr_{\ti k_v}(\eta)-a_1\endpmatrix, 
\pmatrix 1 & \eta^{\sig}\\ \eta & \n_{\ti k_v}(\eta)-a_2\endpmatrix
\right)
\end{equation*}
and so $x$ satisfies (\ref{lcond}), (\ref{licond})
or (\ref{ur*cond}).  
This proves the second implication of (1)--(3).  
\end{proof} 

Our next task is to prove that
the sets defined in (\ref{lcond})--(\ref{*cond})
are covered by the $K_v$-orbits of suitably chosen
standard representatives.
\begin{lem} \label{pplem}
Suppose $p,p'$ are  Eisenstein polynomials whose roots generate 
the same ramified quadratic field over $k_v$.  
Then there exists $\ka\in K_v$ such that $w_p = \ka w_{p'}$.
\end{lem}
\begin{proof}  Let $\al_1,\al_2$ and $\al_1',\al_2'$ be
roots of $p,p'$, respectively.  Since $\al_1,\al_1'$
are both uniformizers of the same field, there exist $c\in \gp_v,\;
d\in \co_v^{\times}$ such that $\al_1=c+d\al_1'$.  Let 
$\ka_{\al,\al'} = n(-c)a(1,d)$.  Then $\ka_{\al,\al'}\in \gl(2)_{\co_v}$
and $h_{\al}= \ka_{\al,\al'}h_{\al'}$.  If $k_v(\al_1)\not=\ti k_v$
then, by \cite{kable-yukie-pbh-I}, (3.18),
\begin{align*}
w_p&=(h_{\al},(\al_2-\al_1)^{-1}h_{\al})w \\
&=(\k_{\al,\al'},d^{-1}\k_{\al,\al'})(h_{\al'},(\al'_2-\al'_1)^{-1}
h_{\al'})w \\
&=(\k_{\al,\al'},d^{-1}\k_{\al,\al'})w_{p'}
\end{align*}
and $(\k_{\al,\al'},d^{-1}\k_{\al,\al'})\in K_v$ since $d\in
\calo_{v}^{\times}$. If $k_v(\al_1)=\ti k_v$ then
\begin{align*}
w_p&=(h_{\al},h_{\al},(\al_2-\al_1)^{-1}h_{\al})w \\
&=(\k_{\al,\al'},\k_{\al,\al'},d^{-1}\k_{\al,\al'})w_{p'}\,.
\end{align*}
Note that we are regarding
$(h_{\al},h_{\al},(\al_2-\al_1)^{-1}h_{\al})$
and $(\ka_{\al,\al'},\ka_{\al,\al'},d^{-1}\ka_{\al,\al'})$ as elements
of $G_{\ti k_v}$  here.  
Since $c,d\in k_v$, $(\ka_{\al,\al'},\ka_{\al,\al'},
d^{-1}\ka_{\al,\al'})$
is an element of $G_{k_v}$ regarded as embedded in $G_{\ti k_v}$.
Therefore
$(\ka_{\al,\al'},\ka_{\al,\al'},d^{-1}\ka_{\al,\al'})\in K_v$ 
in this case also.  
\end{proof} 

If $\ell_1\neq\ti\ell$, we choose Eisenstein polynomials
$p_{\ell_1,j}$, for $j=1,\dots, N_{\ell_1}$, so that
$\{k_v(w_{p_{\ell_1,j}})\}$ is a complete set of representatives for
the classes in $\gX_{\ell_1}$. Similarly, we choose Eisenstein
polynomials $p_{\ti\ell,i,j}$, for $j=1,\dots, N_{\ti\ell}(i)$, so
that $\{k_v(w_{p_{\ti\ell,i,j}})\}$ is a complete set of
representatives for the classes in $\gX_{\ti\ell}(i)$ and an
Eisenstein polynomial $p_{\ti\ell}^{\text{ur}}$ so that
$\gX_{\ti\ell}^{\text{ur}}$ is the singleton containing the class of
$k_v(w_{p_{\ti\ell}^{\text{ur}}})$.  In order to simplify the notation, 
we write $w_{\ell_1,j}$ in place of
$w_{p_{\ell_1,j}}$, $w_{\ti\ell,i,j}$ in place of
$w_{p_{\ti\ell,i,j}}$ and $w_{\ti\ell}^{\text{ur}}$ in place of
$w_{p_{\ti\ell}^{\text{ur}}}$.
\begin{prop}\label{korbit}
\begin{itemize}
\item[(1)] If $x\in {\mathcal D}_{\ell_1}$ then
$x\in \cup_j K_vw_{{\ell_1,j}}$.
\item[(2)] If $x\in {\mathcal D}_{\ti\ell}(i)$ then 
$x\in \cup_j K_vw_{{\ti\ell,i,j}}$.
\item[(3)] If $x\in {\mathcal D}^{\text{\upshape ur}*}_{\ti\ell}$ then
$x\in  K_vw^{\text{\upshape ur}}_{\ti\ell}\cup K_vw_{\eta}$.
\item[(4)] If $x\in {\mathcal D}^*_{\ti\ell}$ then 
$x\in K_vw_{\eta}$.
\end{itemize}
\end{prop}
\begin{proof} 
As shown in the proof of Proposition \ref{xbelong}, 
we may assume that $x_{10}=0$ and $x_{11}=x_{20}=1$. 
Let $p(z)=z^2+a_1(x)z+a_2(x)$.  Then, by (\ref{a1a2}), 
\begin{equation*}
x = (n(x_{21}^{\sigma}-a_1(x)),1)w_p \in K_vw_p
.\end{equation*}

Consider (1). By Proposition \ref{xbelong}
there exists $j$ such that $k_v(x)=k_v(w_{{\ell_1,j}})$.  
By Lemma \ref{pplem}, $w_p\in K_vw_{{\ell_1,j}}$ and
so $x\in K_vw_{{\ell_1,j}}$. Cases (2), (3) and (4) are similar.
\end{proof} 
Next we shall find the volume of the sets 
defined in (\ref{lcond})--(\ref{*cond}) and find the number of 
copies needed to cover the $K_v$-orbits of the
standard representatives.  

\begin{lem}\label{u1u2order}
Suppose $u_1,u_2\in\ti\co_v$, $\ord_{\ti k_v}(u_1) = j$,
and $\ord_{\ti k_v}(u_2)\geq j+1$.   
\begin{itemize}
\item[(1)] If  $j<\ti\del_v$ then $\ord_{k_v}(\n_{\ti k_v/k_v}(u_1))< 
\ord_{k_v}(\tr_{\ti k_v/k_v}(u_2))$. 
\item[(2)] If $j=\ti\del_v$ or $\ti\del_v+1$ then
$\n_{\ti k_v/k_v}(u_1),\tr_{\ti k_v/k_v}(u_2)\in \gp_v^j$.  
\end{itemize}
\end{lem}
\begin{proof}  By Lemma \ref{jorder}, $\ord_{k_v}(\tr_{\ti
k_v/k_v}(u_2))\geq
\lfloor(j+\ti\del_v+1)/2\rfloor\geq(j+\ti\del_v)/2$. If $j<\ti\del_v$
then $j<(j+\ti\del_v)/2$ and, since $\ord_{k_v}(\n_{\ti
k_v/k_v}(u_1))=j$, (1) follows. In (2), it is clear that $\n_{\ti
k_v/k_v}(u_1)\in\gp_v^j$. If $j=\ti\del_v$ then Lemma \ref{jorder}
gives $\ord_{k_v}(\tr_{\ti k_v/k_v}(u_2))\geq\ti\del_v$ and if
$j=\ti\del_v+1$ then it gives $\ord_{k_v}(\tr_{\ti
k_v/k_v}(u_2))\geq\ti\del_v+1$. This completes the verification of
(2).
\end{proof}
\begin{prop}\label{dvolume}
\begin{itemize}
\item[(1)] Suppose $\ell_1\not=\ti\ell$ and 
$\ell_1\leq m_v$ and let $\ell =\operatorname{min}\{\ell_1,\ti\ell\}$.
Then $\vol({\mathcal D}_{\ell_1})= 
q_v^{-\ell_1-\ell-2}(1-q_v^{-1})^4$.  
\item[(2)] Suppose $\ell=\ti\ell\leq m_v$ and $\ell_1=m_v+1$. 
Then $\vol({\mathcal D}_{\ell_1})= 
q_v^{-\ell_1-\ell-2}(1-q_v^{-1})^3$. 
\item[(3)] Suppose $\ti\ell\leq m_v$.  Then 
$\vol({\mathcal D}_{\ti\ell}(\ti\ell))= 
q_v^{-2\ti\ell-2}(1-q_v^{-1})^3(1-2q_v^{-1})$.  
\item[(4)] Suppose $\ti\ell<i<\ti\del_v$ or $\ti\ell=m_v+1$.  
Then $\vol({\mathcal D}_{\ti\ell}(i))= 
q_v^{-2i-2}(1-q_v^{-1})^4$.   
\item[(5)] $\vol({\mathcal D}^{\text{\upshape ur}*}_{\ti\ell})=
q_v^{-2\ti\del_v-1}(1-q_v^{-1})^3$.     
\item[(6)] $\vol({\mathcal D}^*_{\ti\ell}) = q_v^{-8(\ti\del_v+1)}$.  
\end{itemize}
\end{prop}
\begin{proof}  Part (6) is obvious.  

Consider (1) and (2). Suppose $\ell_1<\ti\ell$ and 
$x\in V_{\co_v}$ satisfies the condition (\ref{lcond})
except possibly for the last condition.  Then 
$a_1(x)\equiv 
\tr_{\ti k_v/k_v}(x_{11}x_{21}^{\sig})- x_{12}x_{20}
\; (\gp_v^{\ell+1})$ by (\ref{aform}).    
Since $\gp_v^{\ti\ell}\sub \gp_v^{\ell+1}$, 
$\tr_{\ti k_v/k_v}(x_{11}x_{21}^{\sig}) 
\equiv 0 (\gp^{\ell+1})$, by Lemma \ref{tracelemma},  
and $\ord_{k_v}(x_{12}x_{20})=\ell$.
Thus the last condition of (\ref{lcond}) is automatically 
satisfied.  The volumes of the sets of $x_{10},x_{11},x_{12},
x_{20},x_{21},x_{22}$ satisfying
condition (\ref{lcond}) are 
$1, 1-q_v^{-1},q_v^{-\ell}(1-q_v^{-1}),
1-q_v^{-1},q_v^{-1}(1-q_v^{-1}),q_v^{-\ell-1}$,
respectively.   Therefore 
\begin{equation*}
\vol({\mathcal D}_{\ell_1}) 
= q_v^{-2\ell-2}(1-q_v^{-1})^4
= q_v^{-\ell_1-\ell-2}(1-q_v^{-1})^4
.\end{equation*}

Now suppose that $\ell_1>\ti\ell$ and again assume that $x\in
V_{\calo_v}$ satisfies the conditions of (\ref{lcond}) except possibly
for the last. We have $\ti\ell\leq m_v$ and so, by Lemma
\ref{tracelemma}, $\ord_{k_v}(\tr_{\ti
k_v/k_v}(x_{11}x_{21}^{\sigma}))=\ti\ell$. Since $\ell=\ti\ell$,
$\ord_{k_v}(x_{10}x_{22})\geq\ti\ell+1$ and so 
\begin{equation*}
\ord_{k_v}(\tr_{\ti k_v/k_v}(x_{11}x_{21}^{\sigma})-x_{10}x_{22})
=\ti\ell\,.
\end{equation*}
If $\ell_1\leq
m_v$ then it follows from this and (\ref{aform}) that $x$ satisfies
the last condition of (\ref{lcond}) if and only if
\begin{equation}\label{x12cong}
x_{12}\equiv -x_{20}^{-1}(\tr_{\ti k_v/k_v}(x_{11}x_{21}^{\sigma})-
x_{10}x_{22})\; (\gp_v^{\ell_1})
\end{equation}
but the corresponding congruence with $\ell_1+1$ in place of $\ell_1$
is false. With the other variables fixed, the volume of the set of
$x_{12}$ satisfying (\ref{x12cong}) is $q_v^{-\ell_1}$ and hence the
volume of the set of allowable $x_{12}$ is
$q_v^{-\ell_1}-q_v^{-(\ell_1+1)}=q_v^{-\ell_1}(1-q_v^{-1})$. This
gives
\begin{align*}
\vol(\mathcal{D}_{\ell_1})&=
(1-q_v^{-1})(q_v^{-\ell_1}(1-q_v^{-1}))(1-q_v^{-1})
(q_v^{-1}(1-q_v^{-1}))q_v^{-(\ell+1)} \\
&=q_v^{-\ell_1-\ell-2}(1-q_v^{-1})^4
\end{align*}
in this case. If $\ell_1=m_v+1$ the reasoning is the same, except that
(\ref{x12cong}) is the only condition on $x_{12}$. We thus obtain a
similar formula for $\vol(\mathcal{D}_{\ell_1})$ with one fewer
factors of $(1-q_v^{-1})$. This proves (1) and (2).

Consider (3).  Suppose $\ti\ell\leq m_v$
and $x$ satisfies the conditions of (\ref{licond}) 
for $i=\ti\ell$ except possibly for the last condition.  
Since $x_{10}x_{22}\in \gp_v^{\ti\ell+1}$,
\begin{equation}\label{a1expression}
a_1(x) \equiv \tr_{\ti k_v/k_v}(x_{11}x_{21}^{\sig})
-x_{12}x_{20} \; (\gp_v^{\ti\ell+1})
\end{equation}
by (\ref{aform}).   The order of the first term is 
$\ti\ell$, by Lemma \ref{tracelemma},
and the order of $x_{12}$ is $\ti\ell$.   
So, when $x_{11},x_{21},x_{12}$
are fixed, for $a_1(x)$ to be of order $\ti\ell$, 
$x_{20}$ has to be a unit which is not congruent to 
$\tr_{\ti k_v/k_v}(x_{11}x_{21}^{\sig})x_{12}^{-1}$ modulo $\gp_v$.
Therefore
\begin{equation*}
\begin{aligned}
\vol({\mathcal D}_{\ti\ell}(\ti\ell))
& = (1-q_v^{-1}) (q_v^{-\ti\ell} (1-q_v^{-1}))
(1-2q_v^{-1}) (q_v^{-1}(1-q_v^{-1})) q_v^{-\ti\ell-1} \\
& = q_v^{-2\ti\ell-2}(1-q_v^{-1})^3(1-2q_v^{-1})
.\end{aligned}
\end{equation*}

Consider (4). Suppose that $x\in V_{\calo_v}$ satisfies all the
conditions of (\ref{licond}) except possibly for the last. We shall
show that the last condition follows automatically. First suppose that
$\ti\ell\leq m_v$ and $\ti\ell< i<\ti\del_v$. Then
$\ord_{k_v}(\tr_{\ti k_v/k_v}(x_{11}x_{21}^{\sigma}))=\ti\ell$, by
Lemma \ref{tracelemma}, and
$x_{10}x_{22},x_{12}x_{20}\in\gp_v^{i}\sub\gp_v^{\ti\ell+1}$. 
Thus, by (\ref{aform}), $\ord_{k_v}(a_1(x))=\ti\ell$, as claimed. Now
suppose that $\ti\ell=m_v+1$ and $\ti\ell\leq i<\ti\del_v$. Then
$\ord_{k_v}(\tr_{\ti k_v/k_v}(x_{11}x_{21}^{\sigma}))\geq\ti\ell$, by
Lemma \ref{tracelemma}, and
$x_{10}x_{22},x_{12}x_{20}\in\gp_v^{i}\sub\gp_v^{\ti\ell}$.
Thus, by (\ref{aform}), $\ord_{k_v}(a_1(x))\geq\ti\ell$ and again the
last condition holds.
This implies that 
\begin{equation*}
\begin{aligned}
\vol({\mathcal D}_{\ti\ell}(i))
& = (1-q_v^{-1})(q_v^{-i}(1-q_v^{-1}))(1-q_v^{-1})
(q_v^{-1}(1-q_v^{-1})) q_v^{-i-1} \\
& = q_v^{-2i-2}(1-q_v^{-1})^4
.\end{aligned}
\end{equation*}

Consider (5). Suppose 
$x$ satisfies the condition (\ref{ur*cond}) 
except possibly for the last condition. 
Then $x_{12},x_{22}\in \gp_v^{\ti\del_v}\sub
\gp_v^{\ti\ell+1}$. 
Since $\ord_{k_v}(\tr_{\ti k_v/k_v}(x_{11}x_{21}^{\sig}))
=\ti\ell$ if $\ti\ell\leq m_v$ and $\geq m_v+1$
if $\ti\ell=m_v+1$, so the last condition of (\ref{ur*cond}) 
is always satisfied.  
Therefore, 
\begin{equation*}
\begin{aligned}
\vol({\mathcal D}^{\text{ur}*}_{\ti\ell})
& = (1-q_v^{-1}) q_v^{-\ti\del_v}(1-q_v^{-1})
(q_v^{-1}(1-q_v^{-1})) q_v^{-\ti\del_v} \\
& = q_v^{-2\ti\del_v-1}(1-q_v^{-1})^3
.\end{aligned}
\end{equation*}
This finishes all the cases.
\end{proof}
\begin{prop}\label{copies}
\begin{itemize}
\item[(1)]  Let $\ell_1\not=\ti\ell$ and 
$\ell=\operatorname{min}\{\ell_1,\ti\ell\}$. 
Suppose $g\in K_v$, $x,y\in {\mathcal D}_{\ell_1}$
and $gx=y$.  Then $g\in H(\ell)$.  
\item[(2)]  Let $\ti\ell\leq i<\ti\del_v$.  Suppose 
$g\in K_v,\; x,y\in {\mathcal D}_{\ti\ell}(i)$ and $gx=y$. 
Then $g\in H(i)$. 
\item[(3)]  Suppose 
$g\in K_v,\; x,y\in {\mathcal D}^{\text{\upshape ur}*}_{\ti\ell}$
and $gx=y$. 
Then $g\in H(\ti\del_v-1)$.  
\end{itemize}
\end{prop}
\begin{proof} Suppose $g=(g_1,g_2)$ is as in 
(\ref{gform}).  
Since both $F_x(v)$ and $F_y(v)$ are congruent to 
unit scalar multiples of $v_1^2$ modulo $\gp_v$, 
$g_{221}\equiv 0\;(\gp_v)$.  Since $(1,g_2)\in 
H(j)$ for every $j\geq 0$, 
we may assume that $g_2=1$.  

Since $x_{20},y_{20}\in\co_v^{\times}$, 
$x_{21},y_{21}\in\ti\gp_v$ and $x_{22},y_{22}\in\gp_v$, 
we have $g_{121}\equiv 0 \;(\ti\gp_v)$.  
This implies that $g_{111},g_{122}\in\ti\co_v^{\times}$.  
By (\ref{yform}), 
\begin{equation*}
y_{22} = \n_{\ti k_v/k_v}(g_{121}) x_{20}
+ \tr_{\ti k_v/k_v}(g_{121}g_{122}^{\sig}x_{21})
+ \n_{\ti k_v/k_v}(g_{122}) x_{22}
.\end{equation*}

Consider (1).  Since $\ell<\ti\del_v$, 
if $\ord_{\ti k_v}(g_{121})\leq \ell$ then  
\begin{equation*}
\ord_{k_v}(\tr_{\ti k_v/k_v}(g_{121}g_{122}^{\sig}x_{21}))
> \ord_{k_v}( \n_{\ti k_v/k_v}(g_{121}) x_{20})=\ell
\end{equation*}
by Lemma \ref{u1u2order}.  Since $x_{22},y_{22}\in \gp_v^{\ell+1}$, 
this is a contradiction.  
Since $i,\ti\del_v-1<\ti\del_v$, (2) and (3) are similar. 
\end{proof} 

The next corollary follows easily from  Lemma \ref{hinvariant}
and Proposition \ref{copies}.
\begin{cor}\label{disjoint}
\begin{itemize}
\item[(1)]  Suppose $\ell_1\not=\ti\ell$.
If $g,g'\in K_v$ and 
$g{\mathcal D}_{\ell_1}\cap g'{\mathcal D}_{\ell_1}\not=\emptyset$
then $g{\mathcal D}_{\ell_1} = g'{\mathcal D}_{\ell_1}$.   
\item[(2)]  Suppose $\ti\ell\leq i<\ti\del_v$.
If $g,g'\in K_v$ and 
$g {\mathcal D}_{\ti\ell}(i)\cap g'{\mathcal D}_{\ti\ell}(i)\not
=\emptyset$ then
$g {\mathcal D}_{\ti\ell}(i) = g'{\mathcal D}_{\ti\ell}(i)$.  
\item[(3)]  
If $g,g'\in K_v$ and
$g{\mathcal D}^{\text{\upshape ur}*}_{\ti\ell}\cap 
g'{\mathcal D}^{\text{\upshape ur}*}_{\ti\ell}\not=\emptyset$ then
$g{\mathcal D}^{\text{\upshape ur}*}_{\ti\ell}=
g'{\mathcal D}^{\text{\upshape ur}*}_{\ti\ell}$.  
\end{itemize} 
\end{cor}

We are now ready to calculate $\sum_x \vol(K_v x)$. 
\begin{prop}\label{orbit-volume}
\begin{itemize}
\item[(1)] Suppose $\ell_1\not=\ti\ell$ and $\ell_1\leq m_v$. 
Then $\sum_j \vol(K_v w_{\ell_1,j})= 
q_v^{-\ell_1}(1-q_v^{-1})^2(1-q_v^{-2})^2$.  
\item[(2)] Suppose $\ti\ell\leq m_v$ and $\ell_1=m_v+1$. 
Then $\sum_j \vol(K_v w_{\ell_1,j})= 
q_v^{-\ell_1}(1-q_v^{-1})(1-q_v^{-2})^2$.  
\item[(3)] Suppose $\ti\ell\leq m_v$. 
Then $\sum_j \vol(K_v w_{\ti\ell,\ti\ell,j})= 
q_v^{-\ti\ell}(1-q_v^{-1})(1-2q_v^{-1})(1-q_v^{-2})^2$.  
\item[(4)] Suppose $\ti\ell<i<\ti\del_v$ or $\ti\ell=m_v+1$.  
Then $\sum_j \vol(K_v w_{\ti\ell,i,j})= 
q_v^{-i}(1-q_v^{-1})^2(1-q_v^{-2})^2$.
\item[(5)] $\vol(K_v w^{\text{\upshape ur}}_{\ti\ell}) 
+ \vol(K_v w_{\eta}) 
= q_v^{-\ti\del_v}(1-q_v^{-1})(1-q_v^{-2})^2$.   
\end{itemize}
\end{prop}
\begin{proof} Consider (1). Let 
$\ell=\operatorname{min}\{\ell_1,\ti\ell\}$.
By Propositions \ref{korbit} and \ref{copies}
and  Corollary \ref{disjoint}, 
$\cup_jK_v w_{\ell_1,j}$ is a disjoint union of translates of
${\mathcal D}_{\ell_1}$ and the number of translates is 
$Q(\ell)$. So, by (\ref{hjorder}) and Proposition \ref{dvolume}(1),
\begin{equation*}
\begin{aligned}
\sum_j \vol(K_v w_{\ell_1,j})
& = Q(\ell) \vol({\mathcal D}_{\ell_1})
= q_v^{\ell+2}(1+q_v^{-1})^2 q_v^{-\ell_1-\ell-2}(1-q_v^{-1})^4 \\
& = q_v^{-\ell_1}(1-q_v^{-1})^2(1-q_v^{-2})^2
.\end{aligned}
\end{equation*}
Cases (2)--(5) are similar using   
Proposition \ref{dvolume}(2)--(5).  
\end{proof} 

Our next task is to determine $\vol(K_v w_{\eta})$. Let $\ti p(z)$ be
the polynomial introduced in the second paragraph of this section and
recall, as noted after (\ref{weta}), that
$w_{\eta}=(n(\eta_2),1)w_{\ti p}$. Let
$G_{w_{\eta}\,\calo_v/\gp_v^{\ti\del_v+1}}$ and
$G^{\circ}_{w_{\eta}\,\calo_v/\gp_v^{\ti\del_v+1}}$ be the sets of
$(\calo_v/\gp_v^{\ti\del_v+1})$-valued points of the schemes
$G_{w_{\eta}}$ and $G^{\circ}_{w_{\eta}}$ over $\calo_v$.

\begin{lem} The order of
$G^{\circ}_{w_{\eta}\,\calo_v/\gp_v^{\ti\del_v+1}}$ is
$q_v^{4\ti\del_v+3}(q_v-1)$.
\end{lem}
\begin{proof}
Since $w_{\eta}=(n(\eta_2),1)w_{\ti p}$ and $(n(\eta_2),1)\in K_v$,
$G^{\circ}_{w_{\eta}\,\calo_v/\gp_v^{\ti\del_v+1}}$ and the similarly
defined set $G^{\circ}_{w_{\ti p}\,\calo_v/\gp_v^{\ti\del_v+1}}$ are
conjugate within $G_{\calo_v/\gp_v^{\ti\del_v+1}}$ and so it suffices
to find the order of $G^{\circ}_{w_{\ti
p}\,\calo_v/\gp_v^{\ti\del_v+1}}$. Let
\begin{equation*}
A_{\ti p}(c,d)=\begin{pmatrix}c&-d\\b_1d&c-b_1d\end{pmatrix}\,.
\end{equation*}
It was proved in \cite{kable-yukie-pbh-I}, Lemma 11.2 that
$G^{\circ}_{w_{\ti p}\,\calo_v/\gp_v^{\ti\del_v+1}}$ consists of
elements of the form $(A_{\ti p}(c_1,d_1),A_{\ti p}(c_2,d_2))$ where
$c_1,d_1\in\ti\calo_v/\ti\gp_v^{2(\ti\del_v+1)}$,
$c_2,d_2\in\calo_v/\gp_v^{\ti\del_v+1}$, $\det A_{\ti p}(c_1,d_1)\in
(\ti\calo_v/\ti\gp_v^{2(\ti\del_v+1)})^{\times}$ and $c_2$ and $d_2$ are
related to $c_1$ and $d_1$ by the equation
\begin{equation*}
A_{\ti p}(c_2,d_2)=A_{\ti p}(c_1,d_1)^{-1}A_{\ti
p}(c_1^{\sigma},d_1^{\sigma})^{-1}\,.
\end{equation*}
Note that $\det A_{\ti
p}(c_1,d_1)\in(\ti\calo_v/\ti\gp_v^{2(\ti\del_v+1)})^{\times}$ if and
only if $c_1\in(\ti\calo_v/\ti\gp_v^{2(\ti\del_v+1)})^{\times}$. The
expression for the order follows immediately.
\end{proof}

We denote by $\bar w_{\eta}$ the reduction of $w_{\eta}$ modulo
$\gp_v^{\ti\del_v+1}$ and by $G_{\bar
w_{\eta}\,\calo_v/\gp_v^{\ti\del+1}}$ the stabilizer of $\bar
w_{\eta}$ in $G_{\calo_v/\gp_v^{\ti\del_v+1}}$. Clearly $G_{w_{\eta}\,
\calo_v/\gp_v^{\ti\del_v+1}}$ is a subgroup of $G_{\bar
w_{\eta}\,\calo_v/\gp_v^{\ti\del_v+1}}$.

\begin{lem} We have $[G_{\bar w_{\eta}\,\calo_v/\gp_v^{\ti\del_v+1}}:
G^{\circ}_{w_{\eta}\,\calo_v/\gp_v^{\ti\del_v+1}}]=
2q_v^{\ti\del_v+2\lfloor\ti\del_v/2\rfloor}$.
\end{lem}
\begin{proof}
Let $\bar w_{\ti p}$ denote the reduction of $w_{\ti p}$ modulo
$\gp_v^{\ti\del_v+1}$ and $G_{\bar w_{\ti
p}\,\calo_v/\gp_v^{\ti\del_v+1}}$ the stabilizer of $\bar w_{\ti p}$
in $G_{\calo_v/\gp_v^{\ti\del_v+1}}$. Our first step will be to show
that every right $G^{\circ}_{w_{\ti p}\,\calo_v/\gp_v^{\ti\del_v+1}}$
coset in $G_{\bar w_{\ti p}\,\calo_v/\gp_v^{\ti\del_v+1}}$ has a
representative of the particular form given in (\ref{partiform})
below.

For $x=(x_1,x_2)\in V_{\calo_v/\gp_v^{\ti\del_v+1}}$, let $\xspan(x)$
be the $(\calo_v/\gp_v^{\ti\del_v+1})$-module spanned by $x_1$ and
$x_2$ inside $W_{\calo_v/\gp_v^{\ti\del_v+1}}$. As was stated in
\cite{kable-yukie-pbh-I}, Lemma 11.4, if
$g_1\in\gl(2)_{\ti\calo_v/\ti\gp_v^{2(\ti\del_v+1)}}$ then there
exists $g_2\in\gl(2)_{\calo_v/\gp_v^{\ti\del_v+1}}$ such that
$(g_1,g_2)\in G_{x\,\calo_v/\gp_v^{\ti\del_v+1}}$ if and only if
$\xspan((g_1,1)x)=\xspan(x)$.

Suppose that $g=(g_1,g_2)\in G_{\bar w_{\ti
p}\,\calo_v/\gp_v^{\ti\del_v+1}}$. Since $F_{w_{\ti p}}(v_1,v_2)$
reduces to $v_1^2$ modulo $\gp_v$, $g_{221}\in
(\gp_v/\gp_v^{\ti\del_v+1})$. Using this fact and examining the
second component of $w_{\ti p}$ modulo $\gp_v$, we find that
$g_{121}\in(\ti\gp_v/\ti\gp_v^{2(\ti\del_v+1)})$. This implies that
$g_{111}$ and $g_{122}$ lie in
$(\ti\calo_v/\ti\gp_v^{2(\ti\del_v+1)})^{\times}$. If we put
$c_1=g_{122}$, $d_1=g_{112}$ and $A_{\ti p}(c_2,d_2)=A_{\ti
p}(c_1,d_1)^{-1}A_{\ti p}(c_1^{\sigma},d_1^{\sigma})^{-1}$ then
$(A_{\ti p}(c_1,d_1),A_{\ti p}(c_2,d_2))\in G^{\circ}_{w_{\ti p}\,
\calo_v/\gp_v^{\ti\del_v+1}}$, the $(1,2)$-entry of $A_{\ti
p}(c_1,d_1)g_1$ is $0$ and the $(1,1)$-entry is $\det(g_1)$. Since
$\det(g_1)\in(\ti\calo_v/\ti\gp_v^{2(\ti\del_v+1)})^{\times}$, we may
further multiply on the left by $(A_{\ti p}(\det(g_1)^{-1},0),*)\in
G^{\circ}_{w_{\ti p}\,\calo_v/\gp_v^{\ti\del_v+1}}$ to find a
representative for the right $G^{\circ}_{w_{\ti
p}\,\calo_v/\gp_v^{\ti\del_v+1}}$ coset of $g$ having the form
\begin{equation}\label{partiform}
\left(\begin{pmatrix}1&0\\u&t\end{pmatrix},*\right)
\end{equation}
with $t\in(\ti\calo_v/\ti\gp_v^{2(\ti\del_v+1)})^{\times}$ and
$u\in\ti\calo_v/\ti\gp_v^{2(\ti\del_v+1)}$. It is easy to check that
each coset has exactly one representative in this form.

Since $w_{\eta}=(n(\eta_2),1)w_{\ti p}$ and $(n(\eta_2),1)\in K_v$, it
easily follows that every right
$G^{\circ}_{w_{\eta}\,\calo_v/\gp_v^{\ti\del_v+1}}$ coset in $G_{\bar
w_{\eta}\,\calo_v/\gp_v^{\ti\del_v+1}}$ also has a unique
representative in the form (\ref{partiform}). Our second step will be
to determine when such an element actually lies in $G_{\bar
w_{\eta}\,\calo_v/\gp_v^{\ti\del_v+1}}$. Suppose that $(g_1,g_2)$ is
in the form (\ref{partiform}). Then $(g_1,g_2)\in G_{\bar
w_{\eta}\,\calo_v/\gp_v^{\ti\del_v+1}}$ if and only if
$\xspan((g_1,1)\bar w_{\eta})=\xspan(\bar w_{\eta})$. Computation
gives $(g_1,1)\bar w_{\eta}=(M_1,M_2)$ where
\begin{equation*}
M_1=\begin{pmatrix}0&t^{\sigma}\\t&tu^{\sigma}+t^{\sigma}u
\end{pmatrix}\,\quad
M_2=\begin{pmatrix}1&u^{\sigma}-t^{\sigma}\eta_2\\
u-t\eta_1&uu^{\sigma}-t\eta_1u^{\sigma}
-t^{\sigma}\eta_2u\end{pmatrix}\,.
\end{equation*}
Note that $y=\left(\begin{smallmatrix}y_0&y_1\\y_1^{\sigma}&y_2
\end{smallmatrix}\right)$ is in $\xspan(\bar w_{\eta})$ if and only if
$y_2=0$ and $y_1-y_1^{\sigma}=y_0(\eta_1-\eta_2)$. Thus our element
lies in the stabilizer of $\bar w_{\eta}$ if and only if
\begin{equation}\label{tueq1}
\begin{aligned}
{} &tu^{\sigma}+t^{\sigma}u=t^{\sigma}-t=0 \\
{} &uu^{\sigma}-t\eta_1 u^{\sigma}-t^{\sigma}\eta_2 u=0 \\
{} &(u^{\sigma}-t^{\sigma}\eta_2)-(u-t\eta_1)=\eta_1-\eta_2\,.
\end{aligned}
\end{equation}
Since $t$ must be a unit, the first equation is equivalent to
$t=t^{\sigma}$, $u^{\sigma}=-u$. Using this, the second two equations
become $u^2=t(\eta_1-\eta_2)u$ and $2u=(\eta_1-\eta_2)(t-1)$. Making
use of the second of these, the first is equivalent to
$u(u+(\eta_1-\eta_2))=0$. Thus (\ref{tueq1}) is equivalent to the
system
\begin{equation}\label{tueq2}
t=t^{\sigma}\,,\ u^{\sigma}=-u\,,\ 
u(u+(\eta_1-\eta_2))=0\,,\ 
2u=(\eta_1-\eta_2)(t-1)\,.
\end{equation}
In the analysis of this system it will be convenient to adopt the
usual abuse of notation by which classes in
$\ti\calo_v/\ti\gp_v^{2(\ti\del_v+1)}$ and their representatives in
$\ti\calo_v$ are denoted by the same symbol.

Since $\ord_{\ti k_v}(\eta_1-\eta_2)=\ti\del_v$, the third equation in
(\ref{tueq2}) is equivalent to the condition that either $\ord_{\ti
k_v}(u)\geq\ti\del_v+2$ or $\ord_{\ti
k_v}(u+(\eta_1-\eta_2))\geq\ti\del_v+2$. These two possibilities are
mutually exclusive and it is easy to check that the bijection
$(u,t)\mapsto(u-(\eta_1-\eta_2),t-2)$ carries the set of solutions to
(\ref{tueq2}) satisfying the first inequality onto the set of
solutions satisfying the second. Thus we may assume henceforth that
$\ord_{\ti k_v}(u)\geq\ti\del_v+2$ provided we then double the number
of solutions found. Since $\ord_{\ti k_v}(u)\geq\ti\del_v+2$,
$\ord_{k_v}(u+u^{\sigma})\geq\lfloor(2\ti\del_v+2)/2\rfloor=
\ti\del_v+1$, by Lemma \ref{jorder}, and so $\ord_{\ti
k_v}(u+u^{\sigma})\geq 2\ti\del_v+2$. Thus the second equation in
(\ref{tueq2}) is a consequence of the third and may be deleted from
the system.

Now suppose that $\ti\del_v\leq 2m_v$, so that $\ti\del_v=2\ti\ell$.
Since $\ord_{\ti k_v}(u)\geq\ti\del_v+2$, we may write
$u=(\eta_1-\eta_2)\pi_v\bar u$ with $\bar u\in\ti\calo_v$. The fourth
equation in (\ref{tueq2}) is then equivalent to $t\equiv 1+2\pi_v\bar
u\;(\ti\gp_v^{\ti\del_v+2})$. Thus $t=1+2\pi_v\bar
u+\pi_v^{\ti\ell+1}c$ with $c\in\ti\calo_v$. It follows that
$t-t^{\sigma}=2\pi_v(\bar u-\bar u^{\sigma})+
\pi_v^{\ti\ell}(c-c^{\sigma})$. But $(\bar u-\bar
u^{\sigma}),(c-c^{\sigma})\in\ti\gp_v^{\ti\del_v}$ and so
\begin{equation*}
t-t^{\sigma}\in\ti\gp_v^{2m_v+2+\ti\del_v}+
\ti\gp_v^{2\ti\ell+2+\ti\del_v}\sub
\gp_v^{2\ti\del_v+2}\,.
\end{equation*}
Thus the first equation of (\ref{tueq2}) is a consequence of the third
and fourth. There are thus
$\#(\ti\gp_v^{\ti\del_v+2}/\ti\gp_v^{2\ti\del_v+2})=q_v^{\ti\del_v}$
choices for $u$ and, for each choice of $u$,
$\#(\ti\gp_v^{\ti\del_v+2}/\ti\gp_v^{2\ti\del_v+2})=q_v^{\ti\del_v}$
choices for $t$. This gives $q_v^{2\ti\del_v}$ solutions to
(\ref{tueq2}) with $\ord_{\ti k_v}(u)\geq\ti\del_v+2$. Thus there are
$2q_v^{2\ti\del_v}$ solutions in all in this case.

Finally, we must consider the case where $\ti\del_v=2m_v+1$. We may
assume that the uniformizer, $\pi_v$, has been chosen so that
$\sqrt{\pi_v}\in\ti k_v$. Since $\ord_{\ti k_v}(u)\geq\ti\del_v+2$, we
may write $u=(\eta_1-\eta_2)\pi_v\bar u$ with $\bar u\in\ti\calo_v$.
Again $t\equiv 1+2\pi_v\bar u\;(\ti\gp_v^{\ti\del_v+2})$ and so
$t=1+2\pi_v\bar u+\sqrt{\pi_v}\pi_v^{m_v+1}c$ with $c\in\ti\calo_v$. 
Let us write $\bar u=\bar u_0+\bar u_1\sqrt{\pi_v}+\bar u_2\pi_v$ and
$c=c_0+c_1\sqrt{\pi_v}+c_2\pi_v$ where $\bar u_0,\bar
u_1,c_0,c_1\in\calo_v$ and $\bar u_2,c_2\in\ti\calo_v$. This is
possible since $\ti k_v/k_v$ is ramified. A simple calculation gives
\begin{equation*}
t-t^{\sigma}=4\bar
u_1\sqrt{\pi_v}\pi_v+2c_0\sqrt{\pi_v}\pi_v^{m_v+1}-
2\pi_v^2(\bar u_2^{\sigma}-\bar u_2)+\sqrt{\pi_v}\pi_v^{m_v+2}
(c_2^{\sigma}+c_2)\,.
\end{equation*}
Now $\bar u^{\sigma}-\bar u_2\in\ti\gp_v^{\ti\del_v}$ and
$c_2^{\sigma}+c_2=(c_2^{\sigma}-c_2)+2c_2\in\ti\gp_v^{2m_v}$ and so
the last two terms lie in
$\ti\gp_v^{2\ti\del_v+3}\sub\ti\gp_v^{2\ti\del_v+2}$. Thus
\begin{equation*}
t-t^{\sigma}\equiv 4\bar u_1\sqrt{\pi_v}\pi_v+
2c_0\sqrt{\pi_v}\pi_v^{m_v+1}\;
(\ti\gp_v^{2\ti\del_v+2})
\end{equation*}
and so $t^{\sigma}\equiv t\;(\ti\gp_v^{2\ti\del_v+2})$ if and only if
$\bar u_1\equiv -(\pi_v^{m_v}/2)c_0\;(\gp_v)$. 
Since $u,t$ are determined modulo $\ti\gp_v^{2\ti\del_v+2}$, 
we can regard $\bar u,c$ as elements of $\ti\co_v/\ti\gp_v^{\ti\del_v}$.
There are $q_v^{2\ti\del_v-1}$ pairs $(\bar
u,c)$ satisfying the congruences relating $\bar u_1$ and $c_0$ and
these lead to $q_v^{2\ti\del_v-1}$ pairs $(u,t)$. Thus there are
$2q_v^{2\ti\del_v-1}$ solutions in all.
\end{proof}
\begin{prop}\label{orbit-volume*}
We have
$\vol(K_vw_{\eta})=\tfrac12q_v^{-\ti\del_v-2\lfloor\ti\del_v/2\rfloor}
(1-q_v^{-1})(1-q_v^{-2})^2$.
\end{prop}
\begin{proof}
In light of the previous two lemmas and Proposition \ref{dvolume}(6),
we have
\begin{align*}
\vol(K_v w_{\eta}) &=
q_v^{-8(\ti\del_v+1)}\cdot
\frac{\#G_{\calo_v/\gp_v^{\ti\del_v+1}}}
{\#G_{\bar w_{\eta}\,\calo_v/\gp_v^{\ti\del_v+1}}} \\
&=q_v^{-8(\ti\del_v+1)}\cdot
\frac{(q_v^2-q_v)^2(q_v^2-1)^2(q_v^{2\ti\del_v+1})^4
(q_v^{\ti\del_v})^4}
{2q_v^{\ti\del_v+2\lfloor\ti\del_v/2\rfloor}\cdot
q_v^{4\ti\del_v+3}\cdot(q_v-1)} \\
&=\tfrac12q_v^{-\ti\del_v-2\lfloor\ti\del_v/2\rfloor}
(1-q_v^{-1})(1-q_v^{-2})^2\,.
\end{align*}
\end{proof}
From Propositions \ref{orbit-volume*} and \ref{orbit-volume}(5) we
easily obtain the following.
\begin{cor} \label{orbit-volume-ur}
We have 
\begin{equation*}
\vol(K_v w_{\ti\ell}^{\text{\upshape ur}}) =  q_v^{-\ti\del_v}
(1-\tfrac12q_v^{-{2\lfloor\ti\del_v/2 \rfloor}})
(1-q_v^{-1})(1-q_v^{-2})^2
.\end{equation*}
\end{cor}
This completes the verification of the values of 
$\bar\vep_v(x)$ in Tables 
\ref{table-dyadic-ungrouped} and \ref{table-dyadic-grouped}.

\bibliographystyle{plain}
\bibliography{pbh2}

\end{document}